\documentclass[12pt]{article}

\usepackage{geometry}
\usepackage{graphics}
\usepackage{amssymb}
\usepackage{amsbsy}
\usepackage{amsmath}
\usepackage{amsthm}
\usepackage{extarrows}
\usepackage{dsfont}
\usepackage{float}
\usepackage[enableskew]{youngtab}
\usepackage{ytableau}
\usepackage[all]{xy}
\usepackage{setspace}

\usepackage{microtype}

\usepackage[square]{natbib}

\newtheorem{thm}{Theorem}[section]
\newtheorem{lem}[thm]{Lemma}
\newtheorem{prop}[thm]{Proposition}
\newtheorem{cor}[thm]{Corollary}

\newtheorem{remark}[thm]{Remark}
\newtheorem{dfn}[thm]{Definition}
\newtheorem{example}[thm]{Example}

\definecolor{grey}{gray}{0.41}

\newcommand{\ma}{\mathcal{A}}
\newcommand{\mb}{\mathcal{B}}
\newcommand{\mc}{\mathcal{C}}
\newcommand{\mg}{\mathcal{G}}
\newcommand{\mh}{\mathcal{H}}
\newcommand{\md}{\mathcal{D}}
\newcommand{\ml}{\mathcal{L}}

\newcommand{\none}{n-1}

\newcommand{\llt}{G_{\boldsymbol{ \nu}} (X;q)}
\newcommand{\Mac}{\widetilde H_{\mu / \rho}(X;q,t)}
\newcommand{\bmu}{{\boldsymbol{\mu}}}
\newcommand{\bnu}{{\boldsymbol{\nu}}}
\newcommand{\blambda}{{\boldsymbol{\lambda}}}
\newcommand{\Skew}{\lambda / \rho}
\newcommand{\Skewmu}{\mu / \rho}

\numberwithin{equation}{section}

\begin{document}

\title{Dual Equivalence Graphs Revisited and the Explicit Schur Expansion of a Family of LLT Polynomials}           

\author{Austin Roberts\thanks{Partially supported by DMS-1101017 from the NSF.}\\ 
Department of Mathematics\\ University of Washington\\
 Seattle, WA 98195-4350, USA\\
 }

\date{\today}

\maketitle


\begin{abstract}     
In 2007 Sami Assaf introduced dual equivalence graphs as a method for demonstrating that a quasisymmetric function is Schur positive.  The method involves the creation of a graph whose vertices are weighted by Ira Gessel's fundamental quasisymmetric functions so that the sum of the weights of a connected component is a single Schur function. In this paper, we improve on Assaf's axiomatization of such graphs, giving locally testable criteria that are more easily verified by computers. We further advance the theory of dual equivalence graphs by describing a broader class of graphs that correspond to an explicit Schur expansion in terms of Yamanouchi words. Along the way, we demonstrate several symmetries in the structure of dual equivalence graphs. We then apply these techniques to give explicit Schur expansions for a family of Lascoux-Leclerc-Thibon polynomials. This family properly contains the previously known case of polynomials indexed by two skew shapes, as was described in a 1995 paper by Christophe Carr\'{e} and Bernard Leclerc. As an immediate corollary, we gain an explicit Schur expansion for a family of modified Macdonald polynomials in terms of Yamanouchi words. This family includes all polynomials indexed by shapes with at most three cells in the first row and at most two cells in the second row, providing an extension to the combinatorial description of the two column case described in 2005 by James Haglund, Mark Haiman, and Nick Loehr.
\end{abstract}

\pagebreak

\tableofcontents

\pagebreak


\section{Introduction}\label{intro}

Dual equivalence was developed and applied by Mark Haiman in \citep{Haiman} as an extension of work done by Donald Knuth in \citep{Knuth}. Sami Assaf then introduced the theory of dual equivalence graphs in her Ph.D. dissertation \citep{Assaf07} and subsequent preprint \citep{Assaf}. In these papers, she is able to associate a number of symmetric functions to dual equivalence graphs and each component of a dual equivalence graph to a Schur function, thus demonstrating Schur positivity. More recently, variations of dual equivalence graphs are given for k-Schur functions in \citep{AB} and for the product of a Schubert polynomial with a Schur polynomial in \citep{ABS}.

A key connection between dual equivalence graphs and symmetric functions is the ring of quasisymmetric functions. The quasisymmetric functions were introduced by Ira Gessel in \citep{Gessel} as part of his work on $P$-partitions. Currently there are a number of functions that are easily expressed in terms of Gessel's \emph{fundamental quasisymmetric functions} that are not easily expressed in terms of Schur functions. For example, such an expansion for plethysms is described in \citep{LW}, for Lascoux-Leclerc-Thibon (LLT) polynomials in \citep{HHLRU}, for Macdonald polynomials in \citep{HHL}, and conjecturally for the composition of the nabla operator with an elementary symmetric function in \citep{HHLRU}. An expressed goal of developing the theory of dual equivalence graphs is to create a tool for turning such quasisymmetric expansions into explicit Schur expansions.

Previously, dual equivalence graphs were defined by five \emph{dual equivalence axioms} that are locally testable and one that is not. One of the main results of this paper is to give an equivalent definition using only local conditions, as stated in Theorem~\ref{4plus}. Many graphs, while not satisfying all of these axioms, correspond to Schur positive expansions. In particular, those admitting a morphism onto a dual equivalence graph, as described in Definition \ref{Morphism}, are necessarily Schur positive. In Theorems \ref{KnuthMorph} and \ref{KnuthMorphCompleteness} we give a classification of the set of graphs admitting such a morphism  and obeying the first dual equivalence axiom. In particular, Theorem~\ref{KnuthMorph} gives an explicit Schur expansion for the symmetric functions associated to such graphs in terms of standardized Yamanouchi words.

The paper concludes by applying the above results to LLT polynomials in Theorem~\ref{LLT Decomposition}. LLT polynomials were first introduced in \citep{LLT} as a $q$-analogue to products of Schur functions and were later given a description in terms of tuples of skew tableaux in \citep{HHLRU}. Corollary~\ref{Mac Decomposition} then applies the results of \citep{HHL} to give an explicit combinatorial description for a family of modified Macdonald polynomials. First introduced in \citep{Mac}, Macdonald polynomials are often defined as the set of $q,t$-symmetric functions that satisfy certain orthogonality and triangularity conditions, as is well described in \citep{Macdonald}. Part of the importance of Macdonald polynomials derives from the fact that they specialize to a wide array of well known functions, including Hall-Littlewood polynomials and Jack polynomials (see \citep{Macdonald} for details). In \citep{Hai01}, Mark Haiman used geometric and representation-theoretic techniques to prove that Macdonald polynomials are Schur positive.

In some cases, nice Schur expansions for LLT and Macdonald polynomials are already known. In particular, the set of LLT polynomials indexed by two skew shapes was described in \citep{CL} and \citep{vanL}, and modified Macdonald polynomials indexed by shapes with strictly less than three columns was described in \citep{HHL} (which in turn drew on the earlier work in \citep{CL}, \citep{vanL}). The first combinatorial description of the two column case was given in \citep{Fishel}, but others were subsequently given in \citep{Zab}, \citep{LM}, and \citep{Assaf08}. In addition, an algorithm for finding the Schur expansion of Macdonald polynomials indexed by shapes with at most four cells in the first row and at most two cells in the second row was given in  \citep{ZabrockiSpecialCase}. Finding a combinatorial interpretation for the three column case is still an open problem, though there is a conjectured formula in \citep{HaglundCombinatorialModel}.

This paper is broken into sections as follows. Section \ref{preliminaries} reviews the necessary material on partitions, tableaux, the Robinson-Schensted-Knuth correspondence, and symmetric functions, before giving the necessary background on dual equivalence graphs. Section \ref{structure} is dedicated to further developing the theory of dual equivalence graphs, culminating in a new axiomatization for dual equivalence graphs in Theorem~\ref{4plus}.
Section \ref{LLTs} applies the results of Section~\ref{structure} to LLT polynomials and Macdonald polynomials. The graph structure given to LLT polynomials in \citep{Assaf} is reviewed before Theorem~\ref{LLT graph} classifies the set of LLT polynomials that correspond to dual equivalence graphs. Theorem~\ref{LLT Decomposition} states that said set of LLT polynomials have a Schur expansion indexed by standardized Yamanouchi words. This set strictly contains the set of LLT polynomials indexed by two skew shapes. Corollary~\ref{Mac Decomposition} then gives a Schur expansion for modified Macdonald polynomials indexed by partition shapes with strictly less than four boxes in the first row and strictly less than three boxes in the second row.

\subsection*{Acknowledgements}
\thispagestyle{empty}

We would like to thank David Anderson, Sami Assaf, Jonathan Browder, Nantel Bergeron, William McGovern, and Isabella Novik for inspiring conversations on this topic. We would also like to thank the referees for there thorough feedback and thoughtful review. Particular thanks belong to Sara Billey for her invaluable guidance throughout the process of creating this paper.

\clearpage

\section{Preliminaries}\label{preliminaries}
This section is dedicated to introducing the key notation and definitions that underlie the rest of the paper. Particular attention is given to known results about dual equivalence graphs.

\subsection{Tableaux}\label{Tableaux}
\begin{figure}[hp]
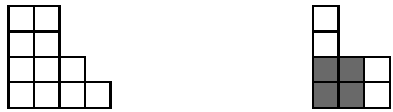

   \begin{center}
    \ytableausetup{smalltableaux}
\ydiagram{2,2,3,4}
\hspace{1in}
\ydiagram[*(white)]{1,1,2+1,2+1}*[*(grey)]{0,0, 2,2}
     \end{center}
  \caption{The diagrams for (4,3,2,2) and  (4,3,1,1)/(2,2) }
 \label{partitions}
\end{figure}

A \emph{partition} $\lambda$ is a weakly decreasing finite  sequence of nonnegative integers $\lambda_1 \geq \ldots \geq \lambda_k  \geq 0$. If $\sum \lambda_i = n$, we say that $\lambda$ is a partition of $n$ and write $\lambda \vdash n$. Partitions are often expressed in terms of diagrams where $\lambda_i$ is the number of boxes, or \emph{cells}, in the $i^{th}$ row, from bottom to top, as in the left diagram of Figure \ref{partitions}. It is sometimes useful to treat a diagram as a subset of the integer Cartesian plane with the bottom left corner of the diagram at the origin. 
 Given a partition $\lambda$, the \emph{conjugate partition} of $\lambda$, denoted $\tilde \lambda$, is defined by $\tilde \lambda_i := |\{j : \lambda_j \geq i\}|$. The diagram of $\tilde \lambda$ is obtained by reflecting the diagram of $\lambda$ over the the main diagonal $x=y$ in the Cartesian plane.

If the diagram of $\rho$ is contained in the diagram of $\lambda$, equivalently $\rho_i \leq \lambda_i$ for all $i$, then we may consider the \emph{skew diagram} $\lambda / \rho$ defined by omitting the boxes of $\rho$ from $\lambda$, as in the right skew diagram of Figure \ref{partitions}. Here, $\Skew$ is referred to as the \emph{shape} of the skew diagram. The number of cells of $\Skew$ is called the \emph{size of} $\Skew$ and is denoted by $|\Skew|$. If we need to distinguish the shape of a partition from a skew shape, we will refer to it as a \emph{straight shape}. We say that $\mu$ is a \emph{subdiagram} of $\Skew$ if some translation of $\mu$ is contained in $\Skew$ when considered as subsets of the Cartesian plane.

\begin{figure}[hp]
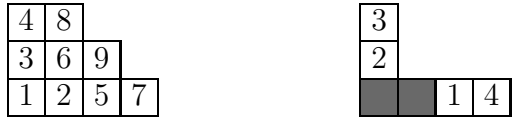

   \begin{center} 
   \ytableausetup{boxsize=.19in}
\ytableaushort{48,369,1257}
\hspace{1in}
\begin{ytableau}
 3  \\
2 \\
*(grey) {}& *(grey) {}&1 &4
\end{ytableau}
   \end{center}
  \caption{On the left, a tableau with content reading word 438162957 and row reading word 483691257. On the right, a skew tableau with content reading word 3214 and row reading word 3214.}
 \label{content}
\end{figure}

A \emph{filling} assigns a positive integer to each cell of a partition or skew shape, usually written inside of the cell. Any filling of $\lambda \vdash n$ that assigns each value in $[n]=\{ 1,\ldots,n\}$ exactly once is termed a \emph{bijective filling}. We will primarily be concerned with \emph{standard Young tableaux}, or \emph{tableaux} for short, which are bijective fillings that are also increasing up columns and across rows from left to right (see Figure \ref{content}). The set of all standard Young tableaux of shape $\lambda$ is denoted SYT($\lambda$). The union of all SYT$(\lambda)$ over $\lambda \vdash n$ is denoted SYT($n$). Similarly, the set of all skew tableaux of shape $\Skew$ is denoted SYT$(\Skew)$.
In general, all tableaux will be assumed to be fillings of straight shapes unless stated otherwise.

The notion of a standard Young tableau extends to fillings of skew shapes, creating \emph{skew tableaux}, as seen in the right side of Figure \ref{content}. The well known process of \emph{jeu de taquin} slides gives a map from any skew tableau to a straight tableau via a sequence of slide operations that move cells either west or south depending on the values in the filling. For more information, see \cite[Part I]{Fulton}, \cite[Ch. 3]{Sagan}, or \cite[Ch. 7]{Stanley2}. Given a (possibly skew) tableau $T$, define sh($T$) to be the shape of the underlying diagram of $T$. 

The \emph{content} of a cell $x$, denoted $c(x)$, is $j-i$, where $j$ is the column of $x$ and $i$ is the row of $x$ in Cartesian coordinates. In other words, each diagonal going southwest to northeast has the same content, with the uppermost diagonal having the smallest content. In a standard Young tableau, the 1-cell is located at the origin of the Cartesian plane, and so has content 0. A connected skew shape having at most one cell of each content is called a \emph{ribbon}.

Define the \emph{content reading word} of a tableau as the word retrieved by reading off each entry from lowest content to highest, moving northeast along each diagonal, as in Figure \ref{content}. We also define the \emph{row reading word} of a tableau by reading across rows from left to right, starting with the top row and working down. The content reading word and row reading word of a standard Young tableau are necessarily permutations.

The \emph{signature} of a word (or permutation in one-line notation) is a string of 1's and -1's, or +'s and $-$'s for short, where there is a + in the $i^{th}$ position if and only if $i$ comes before $i+1$ in the word. If $\sigma_i(w) = -1$, then $i$ is referred to as an \emph{inverse descent of} $w$. Notice that a word is one entry longer than its signature.  We may then define the signature of a tableau $T$, denoted $\sigma(T)$, as the signature of the content reading word of $T$. For example, the signatures of the tableau in Figure \ref{content} are $+--+-+-+$ and $--+$, respectively. Given a fixed tableau, the row reading word and the content reading word always have the same signature (see \cite[Ch.7]{Stanley2} for details).

\subsection{Knuth Equivalence and the R-S-K Correspondence}\label{The RSK Correspondence}

While we assume familiarity with the Robinson-Schensted-Knuth (R-S-K) correspondence and jeu de taquin, we will use this section as a refresher and to set notation. For a full treatment, see \cite[Ch. 2-4]{Fulton}, \cite[Ch. 3]{Sagan}, or \cite[Ch. 7]{Stanley2}.

The R-S-K correspondence  gives a bijection between permutations in $S_n$ and pairs of standard Young tableaux $(P,Q)$, where $P$ and $Q$ have the same shape $\lambda \vdash n$. The first tableau is called the \emph{insertion tableau} and the latter is termed the \emph{recording tableau}. For the duration of this paper, $P\colon S_n\rightarrow \textup{SYT}(n)$ and $Q\colon S_n\rightarrow \textup{SYT}(n)$ will be the functions taking a permutation to its insertion tableau and recording tableau, respectively. These two functions are related by 
\begin{equation}\label{PQ}
Q(w)=P(w^{-1}).
\end{equation}

\noindent A detailed proof of this fact can be found in \cite[Ch. 4.1]{Fulton}. We then write sh($w$) to mean sh$(P(w))$. 

For each tableau $T$ with entries in $[n]$, the set of permutations in $S_n$ sent to $T$ by $P$ is termed a \emph{Knuth equivalence class}. Two words in the same Knuth equivalence class are said to be \emph{Knuth equivalent}. The equivalence relations of Knuth classes are generated by the \emph{fundamental Knuth equivalences}, denoted $K_j$ for $1< j < n$. Each $K_j$ is defined as an involution that fixes all entries of $w\in S_n$ except for those with indices $j-1, j,$ and $j+1$. Its action on these three entries can be written as,
\begin{equation}\label{Knuth Equiv}
\begin{split}
K_j(\ldots xyz \ldots) = (\ldots xyz \ldots), \;\;\;\;\;\;\;
K_j(\ldots zyx \ldots) = (\ldots zyx \ldots),\\
K_j(\ldots yxz \ldots) = (\ldots yzx \ldots), \;\;\;\;\;\;\;
K_j(\ldots xzy \ldots) = (\ldots zxy \ldots),
\end{split}
\end{equation}
where $x < y < z$. In words, if the $j-1,j$, and $j+1$ entries are not strictly increasing or strictly decreasing, then switch the location of the two extreme values.

A number of important constructions yield words from the same Knuth class. Given a tableau $T$ with row reading word $w$, then $P(w)=T$. The same can be shown to be true for the content reading word of $T$, demonstrating that row and content reading words are Knuth equivalent. In fact, the row and content reading words of a skew tableau are also Knuth equivalent. Further, the row reading words (as well as content reading words) of two skew tableaux related by a sequence of jeu de taquin slides are Knuth equivalent. In particular, if $v$ and $w$ are the row reading words of skew tableaux that are related by jeu de taquin, then $\textup{sh}(v)=\textup{sh}(w)$. It also follows that the row reading words of distinct tableaux (on straight shapes) are in different Knuth classes and that there is exactly one such word per class (see \cite[Ch. 2.1]{Fulton} for the details of the proof).

Next, we comment on the relationship between sh$(w)$ and subwords of $w$, as is well presented in \cite[Ch. 3]{Fulton}. If sh$(w)=\lambda$, then the longest increasing subword of $w$ has length $\lambda_1$, and the longest decreasing subword of $w$ has length $\tilde \lambda_1$. For instance, if $w= 15342$, then sh$(w)=(3,1,1)$, the longest increasing subword is 134, and the longest decreasing words are 532 and 542. In particular, if two words are Knuth equivalent, then both of their longest increasing subwords have the same length. 
Furthermore, if $w$ and $v$ are Knuth equivalent words in $S_n$, we may consider the restrictions of $w$ and $v$
 to the consecutive values in some set $S =\{a, a+1, \ldots, b\}$, where $1\leq a<b \leq n$. Call these two subwords $w_S$ and $v_S$, respectively. Then $w_S$ and $v_S$ are Knuth equivalent, and so the longest increasing subwords of $w_S$ and $v_S$ both have the same length. The proof of this last fact can be found in \cite[Lem. 3]{Fulton}.

Lastly, we define a particularly nice Knuth class. Let $U_\lambda$ denote the tableau of shape $\lambda \vdash n$ formed by filling cells with values 1 through $n$ row by row from bottom to top. Define SYam$(\lambda)$ to be the set of $w \in S_n$ such that $P(w)=U_\lambda$. There is, however, a more direct way of deriving this set. A \emph{Yamanouchi word} has entries in the positive integers such that when read backwards there are always more 1's than 2's, more 2's than 3's, and more $i$'s than $i+1$'s. For instance, 25432431121 is a Yamanouchi word, but 231321 is not. The set Yam($\lambda$) consists of all Yamanouchi words where 1 occurs $\lambda_1$ times, 2 occurs $\lambda_2$ times, and so on. We may \emph{standardize} a word in Yam($\lambda)$ by replacing all of the 1's with $1, \ldots, \lambda_1$ in increasing reading order, all of the 2's with $\lambda_1 +1, \ldots, \lambda_2$ in reading order, et cetera. We call the resulting words \emph{standardized Yamanouchi words}. It is a simple exercise to verify that the set of standardized Yamanouchi subwords derived from Yam($\lambda$) is precisely SYam($\lambda)$.

\subsection{Symmetric Functions}\label{Symmetric Functions}

The ring of symmetric functions has several well-known bases with ties to tableaux, as is well laid out in \cite[Ch. 7]{Stanley2}, \cite[Part I]{Fulton}, or \cite[Ch. 4]{Sagan}. Of primary importance is the basis of Schur functions, denoted $\{s_\lambda\}$. We will take the unorthodox approach of defining these functions using a result of Ira Gessel. While less immediately intuitive than standard approaches, this definition contains the only properties that we need. First, a preliminary definition:

\begin{dfn}
\textup{
Given any signature $\sigma \in \{\pm 1\}^{n-1}$, define the \emph{fundamental quasisymmetric function $F_\sigma (X)$} $\in \mathds{Z}[x_1, x_2, \ldots]$ by}
\[
F_\sigma(X):= \sum_{{i_1 \leq \ldots \leq i_n \atop i_j = i_{j+1} \Rightarrow \sigma_j = +1}} x_{i_1}\cdots x_{i_n}.
\]
\end{dfn}

The set of fundamental quasisymmetric functions of degree $n$ forms a homogeneous basis for the vector space of degree $n$ quasisymmetric functions. The \emph{ring of quasisymmetric functions} is created by allowing formal multiplication as power series. The extent that we need this ring to motivate our results is limited to a few facts. The first is the promised definition of Schur functions.

\begin{dfn}\label{Ges}
\textup{{\bf \citep{Gessel}} 
\, Given any skew shape $\Skew$, define 
}
\begin{equation}
s_{\Skew}(X):= \sum_{T \in \textup{SYT}(\Skew)} F_{\sigma(T)}(X),
\end{equation}
\textup{where $s_{\Skew}$ is termed a Schur function if ${\Skew}$ is a straight shape and a skew Schur function in general.}
\end{dfn} 

While it is not obvious from this definition that Schur functions are symmetric or that the Schur functions indexed by straight shapes form a basis for the ring of symmetric functions, what we have gained from this definition is a clear connection to the signatures of tableaux. Further, the quasisymmetric definition of Schur functions is always a finite sum, unlike the more common sum over all semistandard tableaux of a given shape.

The important Lascoux-Leclerc-Thibon (LLT) polynomials and Macdonald polynomials (as introduced in \citep{LLT} and \citep{Mac}, respectively) may also be expressed using the sum of fundamental quasisymmetric polynomials.  
We now present these combinatorial definitions, as they will be needed in Section \ref{LLTs}. The LLT polynomials, denoted $\llt$, were originally described in terms of ribbon tableaux in \citep{LLT}. We will instead use the equivalent definition given in \cite[Cor.~5.2.4]{HHLRU}, which defines $\llt$ by using a $k$-tuple of skew shapes $\bnu$. 

\begin{figure}[h]
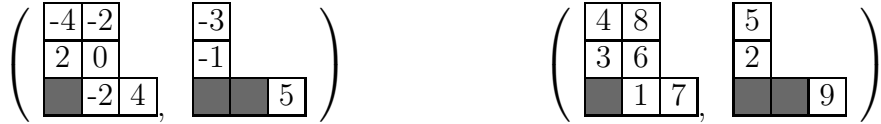

  \[ \ytableausetup{textmode, aligntableaux=bottom}
\left( \begin{array}{cc}  \ytableaushort{{-4}{-2},{2}{0},{*(grey){}}{-2}4} , &
\begin{ytableau}
 {-3}  \\
{-1} \\
*(grey) {}& *(grey) {}&5
\end{ytableau} \end{array} \right)
\hspace{1in}
\left( \begin{array}{cc} \ytableaushort{48,36,{*(grey){}}17}, &
\begin{ytableau}
 5  \\
2 \\
*(grey) {}& *(grey) {}&9
\end{ytableau} \end{array} \right)
\]
  \caption{On the left, the shifted contents of a pair of skew diagrams. On the right, a standard filling of the same tuple with shifted content word 453826179 and signature $---+++-+$.}
\label{shifted content}
\end{figure}

Given a $k$-tuple of skew shapes $\bnu =(\nu^{(0)},\ldots,\nu^{(k-1)} )$, we write $|\bnu|=n$ if $\sum_{i=0}^{k-1}|\nu^{(i)}|=n$. A standard filling ${\bf T} = (T^{(0)},\ldots,T^{(k-1)})$ of $\bnu$ is a bijective filling of the diagram of $\bnu$ with entries in $[n]$ such that for all $0\leq i<k$, each $T^{(i)}$ is strictly increasing up columns and across rows from left to right. Denote the set of standard fillings of $\bnu$ as SYT$(\bnu)$.
Define the \emph{shifted content} of a cell $x$ in $\nu^{(i)}$ as,
\begin{equation}
\tilde c(x) = k\cdot c(x) +i,
\end{equation}

\noindent where $c(x)$ is the content of $x$ in $\nu^{(i)}$. The \emph{shifted content word} of ${\bf T}$ is defined as the word retrieved from reading off the values in the cells from lowest shifted content to highest, reading northeast along diagonals of constant shifted content. We may then define $\sigma({\bf T})$ as the signature of the shifted content word of {\bf T}. For an example, see Figure \ref{shifted content}. 

Letting {\bf T}($x$) denote the entry in cell $x$, the set of $k$-\emph{inversions of} {\bf T} is 
\begin{equation}
\textup{Inv}_k({\bf T}) := \{(x,y) \; |\; k >\tilde c(y) - \tilde c(x) > 0 \, \textup{ and } {\bf T}(x) > {\bf T}(y)\}.
\end{equation}
\noindent The $k$-\emph{inversion number of} {\bf T} is defined as
\begin{equation}
\textup{inv}_k({\bf T}) := |\textup{Inv}_k({\bf T})|.
\end{equation}
If $w$ is the shifted content word of ${\bf T}\in \textup{SYT}(\bnu)$, and $\bnu$ is a $k$-tuple, then the \emph{$\bnu$-inversion number of} $w$ is defined as
\begin{equation}
\textup{inv}_\bnu(w) := \textup{inv}_k({\bf T}).
\end{equation}
\noindent As an example, let {\bf T} be as in Figure \ref{shifted content}. Denoting cells with their values in {\bf T}, Inv$_2({\bf T})$ is comprised of the pairs (5,3), (3,2), and (8,2). Hence,  inv$_2({\bf T})$=inv$_\bnu(453826179)=3$.

 Now define the set of LLT polynomials by
\begin{equation}\label{dfn LLT poly}
\llt := \sum_{{\bf T} \in \textup{SYT}(\bnu)} q^{\textup{inv}_k({\bf T})} F_{\sigma({\bf T})}(X). 
\end{equation}
Though LLT polynomials are known to be symmetric, with proofs in \cite[Thm. 6.1]{LLT} and \cite[Theorem~3.3]{HHL}, it is still challenging to expand them in terms of Schur functions. A partial solution to this problem is given in Section \ref{LLTs}.

We now move on to the definition of the modified Macdonald polynomials $\Mac$. We will use \cite[Theorem~2.2]{HHL} to give a strictly combinatorial definition. To do this, we will first need to define several functions.

\begin{figure}[h]
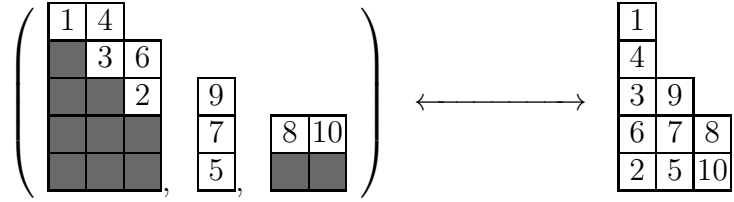

\[
\left( \begin{array}{ccc} \ytableaushort{14,{*(grey){}}36,{*(grey){}}{*(grey){}}2,
{*(grey){}}{*(grey){}}{*(grey){}},{*(grey){}}{*(grey){}}{*(grey){}}}, &
\ytableaushort{9,7,5}, &
\ytableaushort{8{10},{*(grey){}}{*(grey){}}}
\end{array} \right)
\;\; \xleftrightarrow[\hspace{2cm}] \;\;
\begin{array}{c} \ytableaushort{1,4,39,678,25{10}} \end{array}
\]
\caption{An example of the bijection between standard fillings of shapes in TR($\Skewmu$) and bijective fillings of $\Skewmu$.}
\label{TR}
\end{figure}

Given any skew shape $\Skewmu$ with each cell represented by a pair $(i,j)$ in Cartesian coordinates, let TR($\Skewmu$) be the set of tuples of ribbons $\bnu=(\nu^{(0)}, \ldots, \nu^{(k-1)})$, such that $\nu^{(i)}$ has a cell with content $j$ if and only if $(i,-j)$ is a cell in $\Skewmu$. There is then a bijection between standard fillings of shapes in TR($\Skewmu$) and bijective fillings of $\Skewmu$ given by turning each ribbon into a column of $\Skewmu$ as demonstrated in Figure \ref{TR}.

%
Consider any $k$-tuple of ribbon shapes $\bnu =(\nu^{(0)}, \ldots, \nu^{(k-1)})$ with some cell $x\in \nu^{(i)}.$ Define the \emph{arm of $x$} and the \emph{leg of $x$}, denoted $a(x)$ and $l(x)$ respectively, by
%
\begin{equation}
a(x) := |\{ \nu^{(j)}: \, j> i \textup{ and there exists some } y\in \nu^{(j)} \textup{ such that } c(x)=c(y)\}|.
\end{equation}
\begin{equation}
l(x) := |\{ y:\, y\in \nu^{(i)} \textup{ and $c(y) < c(x)$} \}|.
\end{equation} 
 
\noindent Here, $c(x)$ and $c(y)$ always refer to the content within the skew tableaux containing $x$ and $y$, respectively. As an example, if $x$ is the cell containing a six in the left diagram of Figure \ref{TR}, then $a(x)=2$ and $l(x)=3$.

We define a \emph{descent} of $\bnu$ to be any cell $x$ in some $\nu^{(i)}$ that has a cell directly below it and define the descent set of $\bnu$ as 
\begin{equation}
\textup{Des}(\bnu):=\{x \in \bnu:  x  \ \textup{is a descent of }\bnu\}.
\end{equation} 
\noindent For example, the descent set of the tuple of ribbon tableaux in Figure \ref{TR} is the set of cells with values in $\{4, 6, 7, 9\}$. Given a standard filling {\bf T} of $\bnu$, our final three statistics can then be defined as
\begin{equation}
a(\bnu) := \sum_{x\in \textup{Des}(\bnu)} a(x),
\end{equation}
\begin{equation}
\textup{inv}({\bf T}) := \textup{inv}_k({\bf T}) - a(\bnu),
\end{equation}
\begin{equation}
\textup{maj}({\bf T}) :=\textup{maj}(\bnu) :=\sum_{x\in \textup{Des}(\bnu)} 1+ l(x).
\end{equation}

\noindent Using the left diagram in Figure \ref{TR} as an example again, we have $a(\bnu)=3,$ inv$({\bf T})=4-3=1$, and maj({\bf T}) = 9. A simple proof that inv({\bf T}) is always nonnegative can be found in \cite[Sec. 2]{HHL}.

We are now able to define the modified Macdonald polynomials and show their relationship with LLT polynomials:
\begin{equation}\label{Macdonalds}
\Mac := \sum_{ { \bnu \in \textup{TR}(\Skewmu) \atop {\bf T}\in \textup{SYT}(\bnu)} } 
q^{\textup{inv}({\bf T})}t^{\textup{maj}({\bf T})}F_{\sigma({\bf T})} =
\sum_{{\bf \bnu}\in \textup{TR}(\Skewmu)} q^{-a(\bnu)}t^{\textup{maj}(\bnu)}\llt .
\end{equation}

\noindent By using this definition, results about LLT polynomials can be easily translated into results about Macdonald polynomials. 

Lastly, we will have use for the following symmetry of modified Macdonald polynomials. It follows from results in \citep{Macdonald} (see also \cite[Eq. 2.30]{Haglund}) that 
\begin{equation}\label{Mac conjugate}
\Mac = \widetilde H_{\tilde\mu / \tilde\rho}(X;t,q).
\end{equation}

\subsection{Dual Equivalence Graphs}\label{Dual Equivalence Graphs}

We now provide the necessary definitions and results from \citep{Assaf}. We begin by recalling Mark Haiman's dual to the fundamental Knuth equivalences defined in (\ref{Knuth Equiv}). 

\begin{dfn}\textup{
Given a permutation in $S_n$ expressed in one-line notation, define an \emph{elementary dual equivalence} as an involution $d_i$ that interchanges the values $i-1, i,$ and $i+1$ as
}
\begin{equation}\label{dual equiv}
\begin{split}
d_i(\ldots i-1 \ldots i \ldots i+1 \ldots) = (\ldots i-1 \ldots i \ldots i+1 \ldots),\\
d_i(\ldots i+1 \ldots i \ldots i-1 \ldots) = (\ldots i+1 \ldots i \ldots i-1 \ldots),\\
d_i(\ldots i \ldots i-1 \ldots i+1 \ldots) = (\ldots i+1 \ldots i-1 \ldots i \ldots),\\
d_i(\ldots i-1 \ldots i+1 \ldots i \ldots) = (\ldots i \ldots i+1 \ldots i-1 \ldots).
\end{split}
\end{equation}
\textup{Two words are $dual \; equivalent$ if one may be transformed into the other by successive elementary dual equivalences.
}
\end{dfn}

As an example, 21345 is dual equivalent to 51234 because
$d_4(d_3(d_2( 21345 ))) = d_4(d_3( 31245 )) = d_4( 41235 ) = 51234$. Notice that if $i$ is between $i-1$ and $i+1$, then $d_i$ acts as the identity. It follows immediately from (\ref{Knuth Equiv}) and (\ref{dual equiv}) that $d_i$ is related to $K_i$ by  
\begin{equation}\label{eq kd}
d_i(w)=(K_i(w^{-1}))^{-1}.
\end{equation}
By (\ref{PQ}) and (\ref{eq kd}), 
$
Q(w) = P(w^{-1})= P(K_i(w^{-1}))=Q((K_i(w^{-1}))^{-1})=Q(d_i(w)).
$
Thus, $Q$ is constant on dual equivalence classes.

We may also let $d_i$ act on the entries of a tableau $T$ by applying them to the row reading word of $T$. It is not hard to check that the result is again a tableau of the same shape. The transitivity of this action is described in the following theorem.

\begin{thm}[{\cite[Prop.~2.4]{Haiman}}]\label{Haiman}
Two standard Young tableaux on partition shapes are dual equivalent if and only if they have the same shape.
\end{thm}

\noindent If we rewrite Theorem~\ref{Haiman} in terms of permutations, it states that dual equivalence classes are precisely the set of permutations $w$ satisfying $Q(w)=T$ for some fixed tableau T. 

The same action of $d_i$ on tableaux is defined by using the content reading word instead of the row reading word. To see this, recall that the row reading word of a tableau is Knuth equivalent to the content reading word of the same tableau. Given any $w \in S_n$, it follows from \cite[Lemma~2.3]{Haiman} that for all $1 <  i,j < n$,
\begin{equation}\label{eq Hai}
Q(K_j\circ d_i(w))=Q(K_j(w))=Q(d_i\circ K_j(w)).
\end{equation} 

\noindent Applying (\ref{eq kd}) and (\ref{eq Hai}), yields
\begin{equation*}
\begin{array}{rclclcl}
P(d_i\circ K_j(w)) &=& 
Q((d_i\circ K_j(w))^{-1}) & = & Q(K_i (( K_j(w))^{-1})) &&\\ 
&=&Q(K_i\circ d_j(w^{-1})) & = & Q(d_j\circ K_i(w^{-1}))&& \\
& = & Q(d_j (( d_i(w))^{-1})) & = & Q((K_j\circ d_i(w))^{-1}) &=&
P(K_j\circ d_i(w)).
\end{array}
\end{equation*}

\noindent In particular,
\begin{equation}\label{eq commute}
P(d_i\circ K_j(w))=P(K_j\circ d_i(w))=P(d_i(w)).
\end{equation}

\noindent Thus, the fact that the row reading word and content reading word of a tableau are in the same Knuth class implies that they determine the same action of $d_i$ on a tableau.

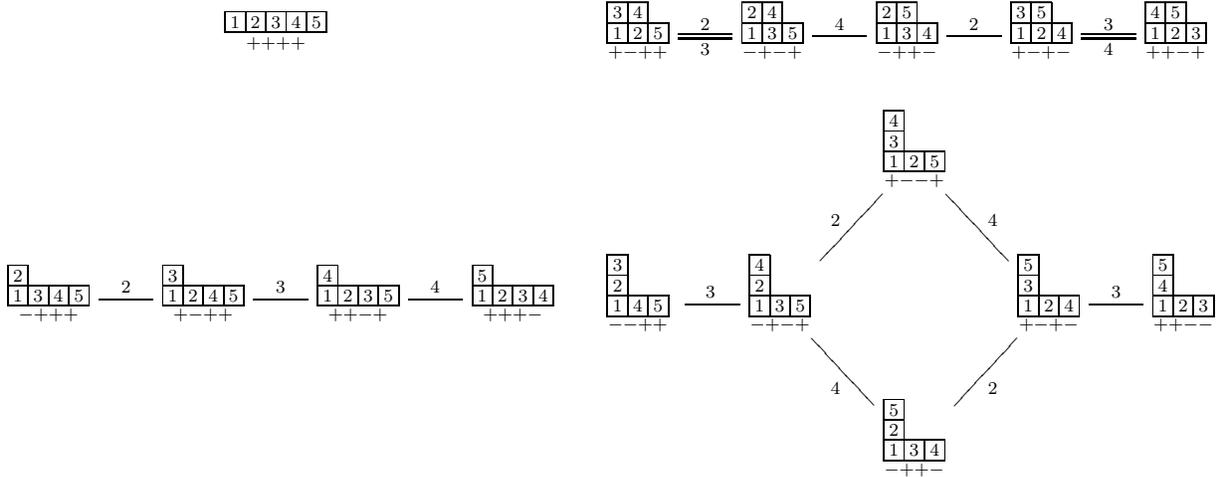
\begin{figure}[H]
\[ \begin{array}{cc}
  \vcenter{\scalebox{.85}{ \vbox{
\xymatrix{ & && {\;\; \Yvcentermath1 \scriptsize {\young(12345) \atop ++++}}}
}}}
& \hspace{-3.7in}
\vcenter{\scalebox{.85}{ \vbox{
	\xymatrix{
		{\Yvcentermath1 \scriptsize {\young(34,125) \atop +-++}} \ar@{=}[r]^{2}_{3} & {\Yvcentermath1 \scriptsize {\young(24,135)}\atop -+-+} \ar@{-}[r]^{4}	& {\Yvcentermath1 \scriptsize {\young(25,134)}\atop -++-} \ar@{-}[r]^2	& { \Yvcentermath1 \scriptsize {\young(35,124)}\atop +-+-} \ar@{=}[r]^3_4 & {\Yvcentermath1 \scriptsize {\young(45,123) \atop ++-+}} 
		} }}}
 \\ & \\
\vcenter{\scalebox{.85}{ \vbox{ \xymatrix{
		{\Yvcentermath1 \scriptsize {\young(2,1345) \atop -+++}} \ar@{-}[r]^{2} 	& {\Yvcentermath1 \scriptsize {\young(3,1245)}\atop +-++} \ar@{-}[r]^{3}	& {\Yvcentermath1 \scriptsize {\young(4,1235)}\atop ++-+} \ar@{-}[r]^4	& { \Yvcentermath1 \scriptsize {{\young(5,1234) \atop +++-}}
} }} } }
&\hspace{-3.7in}
	\vcenter{\scalebox{.85}{ \vbox{ \xymatrix{
		&&{\Yvcentermath1 \scriptsize {\young(4,3,125)}\atop +--+} \ar@{-}[rd]^{4} &&\\
		{\Yvcentermath1 \scriptsize {\young(3,2,145)}\atop --++}\; \ar@{-}[r]^{3}	& {\Yvcentermath1 \scriptsize{\young(4,2,135)}\atop -+-+} \ar@{-}[ur]^{2} \ar@{-}[rd]_{4} && {\Yvcentermath1 \scriptsize{ \young(5,3,124)}\atop +-+-} \ar@{-}[r]^{3} & {\Yvcentermath1 \scriptsize{\young(5,4,123)}\atop ++--}\\
		&&{\Yvcentermath1 \scriptsize {\young(5,2,134)}\atop -++-} \ar@{-}[ur]_{2} &&
		} }} } 
\end{array} \]
\caption{The standard dual equivalence graphs on partitions of 5 up to conjugation.}
 \label{5graph}
\end{figure}

By definition, $d_i$ is an involution, and so we define a graph on standard Young tableaux by letting each nontrivial orbit of $d_i$ define an edge colored by $i$. By Theorem~\ref{Haiman}, the graph on SYT($n$) with edges labeled by $1 < i < n$ has connected components with vertices in SYT($\lambda$) for each $\lambda \vdash n$. We may further label each vertex with its signature to create a \emph{standard dual equivalence graph} that we will denote $\mathcal{G}_\lambda$ (see Figure \ref{5graph}). 

Definition \ref{Ges} and Theorem~\ref{Haiman} determine the connection between Schur functions and dual equivalence graphs as highlighted in \cite[Cor.~3.10]{Assaf}. Given any standard dual equivalence graph $\mg_\lambda=(V,\sigma, E)$,
\begin{equation} \label{schurs}
\sum_{v\in V} F_{\sigma(v)} = s_\lambda.
\end{equation}

\noindent Here, $\mathcal{G}_\lambda$ is an example of the following broader class of graphs.

\begin{dfn}
\textup{
An \emph{edge colored graph} consists of the following data:}
\begin{enumerate}\itemsep-.0in
\item a finite vertex set $V$,
\item a collection $E_i$ of unordered pairs of distinct vertices in $V$ for each \\ $i \in \{m+1, \ldots, n-1\}$, where $m$ and $n$ are positive integers.

\leftskip -0.71in 
\textup{\hspace{.25in} A \emph{signed colored graph} is an edge colored graph with the following additional data:}

\leftskip 0in
\item a signature function $\sigma\colon V \rightarrow \{\pm1\}^{N-1}$ for some positive integer $N \geq n$.
\end{enumerate}

\noindent \textup{We denote a signed colored graph by $\mg =(V, \sigma, E_{m+1}\cup \cdots \cup E_{n-1})$ or simply $\mg=(V, \sigma, E)$.
If a signed colored graph has $m=1$, as described above, then it is said to have \emph{type $(n,N)$} and is termed an $(n,N)$-signed colored graph.
}
\end{dfn}

For our purposes, whenever $V$ is a set of permutations or tableaux, it will be assumed that $\sigma$ is the signature function defined in Section \ref{Tableaux}. To be explicit, we will sometimes refer to this definition of the signature function as \emph{given by inverse descents}.

Signed colored graphs of different types may often be related by restricting some of the data. For example, if $\mg$ is an $(n,N)$-signed colored graph, $M\leq N$, and $m\leq n$, then the \emph{$(m,M)$-restriction} of $\mg$ is the result of excluding $E_i$ for $i \geq m$ and projecting each signature onto its first $M-1$ coordinates.
The $(m,M)$-component of a vertex $v$ of $\mg$ is the connected component containing $v$ in the $(m,M)$-restriction of $\mg$.

In order to describe which signed colored graphs have the same structure as a standard dual equivalence graph, we first need to define isomorphisms.

\begin{dfn}\label{Morphism}
 \textup{
A map $\phi\colon \mg \rightarrow \mh$  between edge colored graphs $\mg=(V, E_{m+1} \cup \ldots \cup E_{n-1})$ and $\mh=(V^\prime, E^\prime_{m+1} \cup \ldots \cup E^\prime_{n-1})$  is called a $morphism$ if it preserves $i$-edges. That is, $\{v, w\}\in E_i$ implies $\{\phi(v), \phi(w)\}\in E^\prime_i$ for all $v, w \in V$ and all $m< i< n$.\\
\indent A map $\phi\colon \mg \rightarrow \mh$  between signed colored graphs $\mg=(V, \sigma, 
 E_{m+1} \cup \ldots \cup E_{n-1})$ and $\mh=(V^\prime,\sigma^\prime, E^\prime_{m+1} \cup \ldots \cup E^\prime_{n-1})$ is called a \emph{morphism} if it is a morphism of edge colored graphs that also preserves signatures. That is, $\sigma^\prime(\phi(v))= \sigma(v)$. \\
\indent In both cases, a morphism is an $isomorphism$ if it admits an inverse morphism.}
\end{dfn}

\noindent Though the term morphism is given two different definitions above, the specific definition should be clear from the context.

The next proposition can be thought of as stating that standard dual equivalence graphs are unique up to isomorphism and have trivial automorphism groups.

\begin{prop}[\citep{Assaf} Proposition 3.11]\label{unique}
If $\phi\colon \mathcal{G}_\lambda \rightarrow \mathcal{G}_\mu$ is an isomorphism of signed colored graphs, then $\lambda=\mu$, and $\phi$ is the identity morphism.
\end{prop}

Notice that in a standard dual equivalence graph, a vertex $v$ is included in an $i$-edge if and only if $\sigma(v)_{i-1}=-\sigma(v)_i$, motivating the following definition.

\begin{dfn}
\textup{Let $\mg = (V,\sigma,E)$ be a signed colored graph. We say that $w \in V$ \emph{admits an $i$-neighbor} if $\sigma(w)_{i-1} = -\sigma(w)_i$.}
\end{dfn}

 Before moving on to an abstract generalization of the structure inherent in any standard dual equivalence graph, recall that a \emph{complete matching} is a simple graph such that every vertex is contained in exactly one edge.

\begin{dfn} \label{axioms}  
\textup{A signed colored graph $\mathcal{G} = (V, \sigma, E_{m+1} \cup \ldots\cup E_{n-1})$ is a \emph{dual equivalence graph}  if the following axioms hold:}
\end{dfn}

\leftskip 0.5in
\parindent -0.5in
(ax1): For $m<i<n$, each $E_i$ is a complete matching on the vertices of $V$ that admit an  $i$-neighbor.

(ax2): If $\{v,w\} \in E_i,$  then $\sigma(v)_i = -\sigma(w)_i, \sigma(v)_{i-1} = -\sigma(w)_{i-1}$, and $\sigma(v)_h = \sigma(w)_h$ for all $h < i-2$ and all $h > i+1$.

(ax3): For $\{v,w\} \in E_i$, if $\sigma_{i-2}$ is defined, then $v$ or $w$ (or both) admits an $(i-1)$-neighbor, and if $\sigma_{i+1}$ is defined, then $v$ or $w$ (or both) admits an $(i+1)$-neighbor.

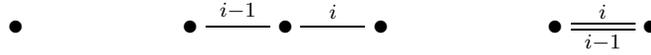
\begin{figure}[H]
  \begin{displaymath}
	\xymatrix{
		\bullet && \bullet \ar@{-}[r]^{i-1} 	& \bullet \ar@{-}[r]^{i}		& \bullet && \bullet \ar@{=}[r]^i_{i-1}	& \bullet
		}
\end{displaymath}  \caption{Allowable $E_{i-1} \cup E_i$  components of Axiom 4}
 \label{4(1)}
\end{figure}

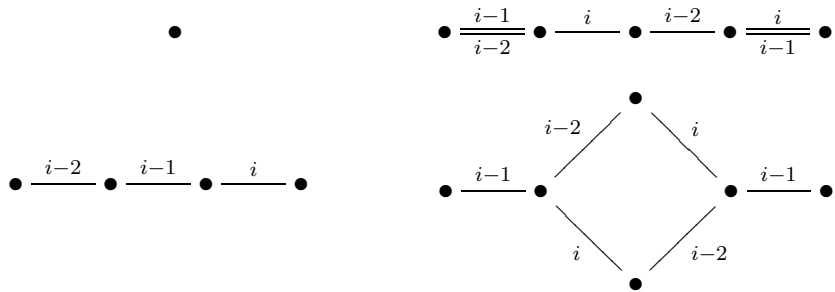
\begin{figure}[H]
\begin{displaymath}
	\xymatrix{
		&&\bullet
		}
	\xymatrix{	
		&&&\bullet \ar@{=}[r]^{i-1}_{i-2}	& \bullet \ar@{-}[r]^i	& \bullet \ar@{-}[r]^{i-2}	& \bullet \ar@{=}[r]^i_{i-1}	&\bullet
		}
\end{displaymath}
\begin{displaymath}
	\xymatrix{
		&&&&\\
		\bullet \ar@{-}[r]^{i-2} 	& \bullet \ar@{-}[r]^{i-1}	& \bullet \ar@{-}[r]^i	& \bullet &
		}
\hspace{.175in}
	\xymatrix{
		&& \bullet \ar@{-}[rd]^{i} &&\\
		\bullet \ar@{-}[r]^{i-1}	& \bullet \ar@{-}[ur]^{i-2} \ar@{-}[rd]_{i} && \bullet \ar@{-}[r]^{i-1} & \bullet\\
		&&\bullet \ar@{-}[ur]_{i-2} &&
		}
\end{displaymath}
\caption{The allowable $E_{i-2} \cup E_{i-1} \cup E_i$  components of Axiom 4.}
 \label{4(2)}
\end{figure}

(ax4): For all $m+1 < i < n$, any component of the edge colored graph $(V, E_{i-2} \cup E_{i-1} \cup E_{i})$ is isomorphic to a component of the restriction of some $\mg_\lambda=(V^\prime, \sigma^\prime, E^\prime)$ to $(V^\prime,E^\prime_{i-2}\cup E^\prime_{i-1} \cup E^\prime_{i})$, where $E_{i-2}$ is omitted if $i= m+2$ (see Figures \ref{4(1)} and \ref{4(2)}).

(ax5): For all $1<i, j< n$ such that $|i-j|>2$, if $\{v,w\} \in E_i$ and $\{w,x\}  \in E_j$, then there exists $y\in V$ such that $\{v,y\} \in E_j$ and $\{x,y\} \in E_i$.

(ax6): For all $m < i < n$, any two vertices of a connected component of $(V, \sigma, E_{m+1} \cup \cdot \cdots \cup E_i)$ may be connected by some path crossing at most one $E_i$ edge.

\leftskip 0.0in \parindent 0.25in

\vspace{.1in}

\noindent A dual equivalence graph that is also an $(n,N)$-signed colored graph is said to have \emph{type $(n,N)$} and is termed an \emph{$(n,N)$-dual equivalence graph}.
%
%
%

\vspace{-.1in}

\begin{remark}\label{axiom remark}
\textup{
 The following are immediate consequences of Definition \ref{axioms}:
\begin{enumerate}\itemsep-.2in
\item A connected component of a dual equivalence graph is also a dual equivalence graph.\\
\item For any $m \leq n$ and $M \leq N$, the $(m,M)$ restriction of an $(n,N)$-dual equivalence graph is an $(m,M)$-dual equivalence graph.\\
\item If a signed colored graph has type $(n,n)$ or if $m+1 < i < n-1$, then Axiom 3 is implied by Axioms 1, 2, and 4 on components of two consecutive colors. In the presence of Axioms 1 and 2, Axiom 3 can be restated in terms of signatures as follows. For $\{v,w\} \in E_i$, if $\sigma(v)_{i-2} =-\sigma(w)_{i-2}$, then $\sigma(v)_{i-2} = -\sigma(v)_{i-1}$ whenever $i>2$, and if $\sigma(v)_{i+1} =-\sigma(w)_{i+1}$, then $\sigma(v)_{i+1} = -\sigma(v)_i$ whenever $\sigma_{i+1}$ is defined. This is the original definition of Axiom 3 used in \textup{\citep{Assaf}}.
\\
\item It is an instructional exercise to check that if Axioms 1, 2, and 6 are obeyed,  then Axiom 4 on 2 consecutive colors implies Axiom 4 on 3 consecutive colors.\\  
\item A signed colored graph $\mg =  (V, \sigma, E_{m+1} \cup \ldots \cup E_{n-1})$ satisfies Axioms 1, 2, 3, and 4, if and only if for any $m< i < n$, each component of $(V, \sigma,  E_{i-2} \cup E_{i-1} \cup E_{i})$ is isomorphic to a component of the restriction of some $\mg_\lambda=(V^\prime, \sigma^\prime, E^\prime)$ to $(V^\prime, \sigma^\prime, E^\prime_{i-2} \cup E^\prime_{i-1} \cup E^\prime_{i})$, where $E^\prime_{i-2}$ or $E^\prime_{i-1}$ is omitted if $i\leq m+2$ or $i=m+1$, respectively. While this fact could be used to shorten the axiomatization, in practice it is often necessary to check Axioms 1, 2, 3, and 4 separately.
\end{enumerate}
}
\end{remark}

The next two theorems link the definition of dual equivalence graphs with that of standard dual equivalence graphs.

\begin{thm}[{\cite[Theorem~3.5]{Assaf}}] \label{Tableaux are DEGs}
For any $\lambda \vdash n$, $\mathcal{G}_\lambda$ is an $(n,n)$-dual equivalence graph.
\end{thm}

\noindent The converse is also true.

\begin{thm}[{\cite[Theorem~3.9]{Assaf}}] \label{main}
Every connected component of an $(n,n)$-dual equivalence graph is isomorphic to a unique $\mg_\lambda.$
\end{thm}

The key to proving Theorem~\ref{main} is building an appropriate morphism. Some of the same techniques will prove useful in this paper, and so we lay them out now.

\begin{dfn} \textup{
Fix any partitions $\lambda \subset \mu$ with $|\lambda| = n$ and $|\mu|= N$ and a skew tableau $A$ of shape $\mu / \lambda$ with entries $n+1, \ldots, N$.  Define the set of standard Young tableaux \emph{augmented  by } $A$, denoted ASYT($\lambda, A$), as the set of $T\in \textup{SYT}(\mu$) such that the restriction of $T$ to $\mu / \lambda$ is $A$. Further, define a signed colored graph $\mg_{\lambda,A}$ on ASYT$(\lambda, A)$ with signature given by inverse descents and edges given by the the nontrivial orbits of $d_i$ for $1< i < n$.
}
\end{dfn}

\begin{remark} \textup{
The graph $\mg_{\lambda,A}$ is isomorphic to an $(n,N)$-component of $\mg_{\mu}$, and every $(n,N)$-component of $\mg_{\mu}$ is isomorphic to some $\mg_{\lambda,A}$, as is clear from the definition of $\mg_{\lambda,A}$. By Part~1 of Remark~\ref{axiom remark} and Theorem~\ref{Tableaux are DEGs}, $\mg_{\lambda,A}$ is therefore a dual equivalence graph  with $(n,n)$-restriction isomorphic to $\mathcal{G}_\lambda$. Applying Theorem~\ref{main}, every $(n,N)$-component of an $(N,N)$-dual equivalence graph is isomorphic to $\mg_{\lambda,A}$ for some $\lambda \vdash n$ and some $A$ such that $|A|=N-n.$ 
}
\end{remark}

With the notion of augmentation it is possible, in some sense, to reverse the process of restriction.

\begin{prop}[{\cite[Lemma~3.13]{Assaf}}]\label{extension} 
Let $\mathcal{G} = (V,\sigma,E)$ be a connected $(n,N)$-dual equivalence graph, and let $\phi$ be a morphism from the $(n, n)$-restriction of $\mathcal{G}$ to $\mathcal{G}_{\lambda}$ for some partition $\lambda$ of $n$. Then $\phi$  extends to an isomorphism $\tilde \phi \colon  \mathcal{G} \rightarrow \mathcal{G}_{\lambda, A}$, where $A$ is a skew tableau such that $|\lambda| + |A| = N$. Furthermore, the position of the cell containing $n+1$ in $A$ is unique.
\end{prop}

Because of the uniqueness statement in Proposition \ref{extension}, we can unambiguously refer to \emph{the unique extension} of a connected $(n,n+1)$-dual equivalence graph to a connected $(n+1,n+1)$-dual equivalence graph. That is, if an $(n,n+1)$-dual equivalence graph is as in Proposition \ref{extension}, then the unique extension is isomorphic to $\mg_\mu$, where $\mu$ is the union of $\lambda$ and the cell of $A$ containing $(n+1)$.

\begin{dfn}
\textup{Let $\mg$ be a signed colored graph of type $(n+1, n+1)$. Two distinct components of the $(n,n+1)$-restriction of $\mg$ that are connected by an $n$-edge in $\mg$ are said to be \emph{neighbors} in $\mg$.}
\end{dfn}

While not explicitly stated in \citep{Assaf}, the following is an immediate consequence of the proof of Theorem~\ref{main}.

\begin{cor}[{\cite[Theorem~3.14]{Assaf}}] \label{complete}
Let $\mg$ be a connected $(n+1,n+1)$-signed colored graph satisfying Axioms 1--5 whose $(n,n+1)$-restriction is a dual equivalence graph. Let $\mc$ be any component of the $(n,n+1)$-restriction of $\mg$, and let the unique extension of $\mc$ be isomorphic to $\mg_\mu$. Then $\mc \cup ( \bigcup {\mb} )$
 is isomorphic to the $(n,n+1)$ restriction of $\mg_\mu$, where the union is over all $\mb$ that are neighbors of $\mc$ in $\mg$.  Furthermore, there exists a morphism $\phi\colon  \mg \rightarrow \mg_\mu$.

\end{cor}

\section{The Structure of Dual Equivalence Graphs}\label{section: structure of DEGs}\label{structure}

The main results of this section are the classification of graphs satisfying Axiom 1 that admit a morphism onto a dual equivalence graph in Theorems \ref{KnuthMorph} and Theorem~\ref{KnuthMorphCompleteness}, the improved axiomatization of dual equivalence graphs given in Theorem~\ref{4plus}, and the more specific criterion for satisfying the dual equivalence axioms given in Corollary~\ref{computer check}. In the process, a number of smaller results about the structure of dual equivalence graphs are highlighted.


\subsection{Symmetries of Dual Equivalence Graphs}

We begin by giving notation for a useful signed colored graph. Let $\mg_n$ denote the $(n,n)$-signed colored graph with vertices indexed by the permutations in $S_n$, signature function given by inverse descents, and $i$-edges given by the nontrivial orbits of $d_i$ for each $1<i<n$.

The following lemma is a natural extension of  \cite[Lemma~2.3]{Haiman}. It lays out a fundamental relationship between the dual equivalence and Knuth equivalence maps defined in (\ref{Knuth Equiv}) and (\ref{dual equiv}).

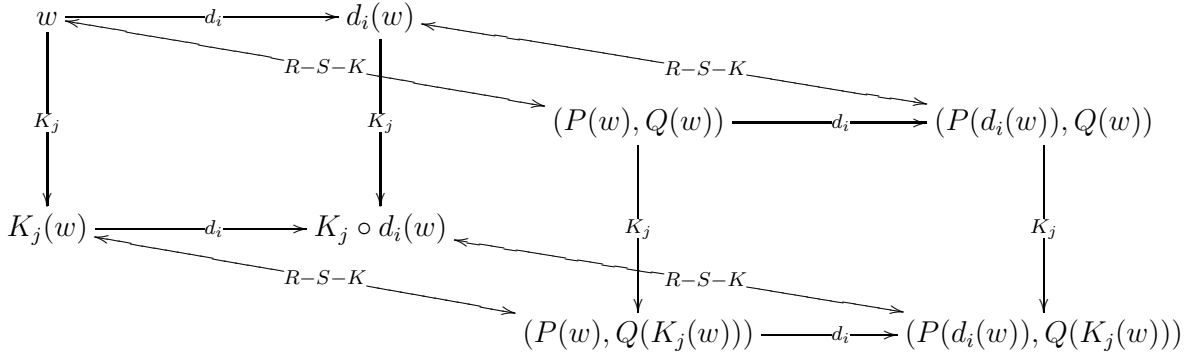
\begin{figure}[H]
   \begin{center} 
   \begin{displaymath}
	\scalebox{.95}{
	\xymatrix{
		w \ar[rrr]|{d_i} \ar[dd]|{K_j} \ar@{<->}[rrrrd]|(.47){R-S-K} &&& d_i(w) \ar[dd]|{K_j} \ar@{<->}[rrrd]|(.49){R-S-K}&&&\\
		&&&&(P(w), Q(w)) \ar[rr]|{d_i} \ar[dd]|{K_j} && (P(d_i(w)), Q(w)) \ar[dd]|{K_j}\\
		K_j(w) \ar[rrr]|{d_i} \ar@{<->}[rrrrd]|(.47){R-S-K}&&& K_j\circ d_i(w) \ar@{<->}[rrrd]|(.49){R-S-K}&&&\\
		&&&&(P(w), Q(K_j(w))) \ar[rr]|{d_i} && (P(d_i(w)), Q(K_j(w)))
		}}
	\end{displaymath}   
   \end{center}
 	 \caption{The commutativity of $d_i$, $K_j$, and the R-S-K correspondence acting on a word $w\in S_n$.}
 \label{Fig Kjdi}
\end{figure}

\begin{lem}\label{Knuth} 
Given any $w\in S_n$ and any $1<i<n$, $1<j< n$, then $K_j\circ d_i(w)=d_i\circ K_j(w)$ and $\sigma(w)=\sigma(K_j(w))$. In particular, $K_j$ defines an automorphism of $\mg_n$.\end{lem}
\noindent $Proof$: 
The fact that $K_j$ preserves inverse descent sets follows from its definition in (\ref{Knuth Equiv}). Thus, $\sigma(w)=\sigma(K_j(w))$. Now we prove that $K_j$ commutes with $d_i$. Recall that the R-S-K correspondence provides a bijection that sends $w\in S_n$ to a pair of tableaux $(P(w),Q(w))$. Consider the effect that applying $d_i$ and $K_j$ to $w$ has on the pair $(P(w),Q(w))$.  By (\ref{eq Hai}) and (\ref{eq commute}), 
\begin{equation}
(P(d_i\circ K_j(w)), Q(d_i\circ K_j(w))) = (P(K_j\circ d_i(w)), Q(K_j\circ d_i(w))).
\end{equation}
Applying the inverse R-S-K correspondence to $S_n$ yields $K_j\circ d_i(w)=d_i\circ K_j(w)$, as shown in Figure \ref{Fig Kjdi}. 

To prove the last part of the lemma, notice that the above argument demonstrates that $K_j$ defines a morphism on $\mg_n$. Because $K_j$ is its own inverse, the morphism must be an isomorphism.~\hfill~$\square$

Knuth equivalences, jeu de taquin, standard skew tableaux, and dual equivalence graphs are all intimately related, as is demonstrated in the next theorem. 

\begin{thm}\label{Skew} 
The function $P$\hspace{0in}$: S_n\rightarrow \textup{SYT}(n)$ induces a surjective morphism from $\mg_n$ to $\cup_{\lambda \vdash n} \mg_\lambda$. This morphism restricts to an isomorphism from any given component of $\mg_n$ to $\mg_\lambda$ for some $\lambda \vdash n$. In particular, if $\Skew$ is a fixed skew shape and $V$ is the set of row reading words of skew tableaux in \textup{SYT}$(\Skew)$, then the restriction of $\mg_n$ to $V$ is a dual equivalence graph.
\end{thm}
\noindent $Proof:$  
We begin by showing that $P$ induces a morphism. As mentioned in Section \ref{The RSK Correspondence}, each $w \in S_n$ is in a Knuth class with the reading word of some $T\in \textup{SYT}(\lambda)$. Thus, there is some sequence of fundamental Knuth equivalences, call it $K_w$, that takes $w$ to the row reading word of $T$. Recall from Sections \ref{Tableaux}, \ref{The RSK Correspondence}, and \ref{Dual Equivalence Graphs} that $P(w) = P(K_w(w)) = T$, $d_i(P(K_w(w))=P(d_i\circ K_w(w))$, and $\sigma(P(w))=\sigma(K_w(w))$. By Lemma~\ref{Knuth}, $\sigma(K_w(w))=\sigma(w)$, so $P$ preserves signatures. Treating $P(w)$ as a vertex in $\mg_\lambda = (V,\sigma, E)$, $P$ takes the vertices in some $i$-edge $\{w, d_i(w)\}$ in $\mg_n$, to 
\begin{align*}
\begin{split}
\{P(w), P(d_i(w))\}  & = \{P(K_w(w)) , P(K_w\circ d_i(w))\} \\  
& = \{P(K_w(w)) , P(d_i\circ K_w(w))\} \\ 
& = \{P(K_w(w)) , d_i(P(K_w(w)))\} \in E_i,
\end{split}
\end{align*}
\noindent again by Lemma~\ref{Knuth}. Therefore, $P$ induces a morphism from $\mg_n$ to $\cup \mg_\lambda$. This morphism is surjective because the R-S-K correspondence guarantees every $T \in \textup{SYT}(n)$ is the image of some $w \in S_n$ under the action of $P$.

To prove the second statement in the theorem, it suffices to restrict the domain to any component $\mc$ in $\mg_n$ and then explicitly create an inverse morphism from $\mg_\lambda$ to $\mc$. Because $Q$ is constant on dual equivalence classes, $Q$ evaluates to some fixed tableau $U$ on all of $\mc$. For any $w \in \mc$, $P(w)=T$ if and only if the inverse R-S-K correspondence sends $(T,U)$ to $w$. The desired inverse morphism is thus given by sending the vertex $T$ to the word corresponding to $(T,U)$ in the R-S-K correspondence. This action takes $P(w)$ to $w$ and $d_i(P(w))=P(d_i(w))$ to $d_i(w)$, so it must preserve edges. The same analysis as the previous paragraph demonstrates that this action preserves signatures, so we have defined the desired  inverse morphism.

For the last statement in the theorem, one can observe from the definition in (\ref{dual equiv}) that SYT$(\Skew)$ is closed under dual equivalence for any fixed shape $\Skew$. Hence, restricting $\mg_n$ to the vertex set $V$ is a restriction to a collection of connected components of $\mg_n$. Thus, $P$ takes each of these components to a standard dual equivalence graph, completing the proof.
\hfill $\square$

\vspace{.1in}

 In light of Theorem $\ref{Skew}$, it makes sense to extend the notation $\mg_\lambda$ to the skew case $\mg_{\Skew}$. That is, $\mg_{\Skew}$ is the dual equivalence graph with vertices in SYT($\Skew$), edges given by nontrivial orbits of $d_i$, and signatures given by inverse descents. As with standard dual equivalence graphs, the actions of $d_i$ and $\sigma$ are defined via the row reading words of skew tableaux. 

\begin{remark}\label{skew via words}
\textup{In the definition of $\mg_{\Skew}$, both $\sigma$ and $E_i$ are defined via row reading words. Thus, $\mg_{\Skew}$ is isomorphic to the signed colored graph induced by sending each vertex to its row reading word. As with standard dual equivalence graphs, both $\sigma$ and $E_i$ may be equivalently defined via content reading words. Thus, $\mg_{\Skew}$ is also isomorphic to the signed colored graph induced by sending each vertex to its content reading word.} 
\end{remark}

Theorem $\ref{Skew}$ leads to a number of simple corollaries. The first is a well-known fact (see \citep{Fulton}, \citep{Sagan}, or \citep{Stanley2}), while the rest help to illuminate the structure of dual equivalence graphs.

\begin{cor}[Littlewood-Richardson Rule]
The skew Schur functions $s_{\nu/\lambda}$ are Schur positive. Moreover, for all $\Skew$,
\[
s_{\Skew} \;
= \sum_{\mu \, \vdash |\lambda|-|\rho|} c_{\rho,\mu}^\lambda \, s_\mu,
\]
where
$c_{\rho, \mu}^\lambda = |\{ w \in \textup{SYam}(\mu): w \textup{ is the row reading word of a skew tableau in SYT}(\Skew)\}|$.

\end{cor}
\noindent $Proof:$ By (\ref{schurs}) and Theorem~\ref{Skew}, we can interpret $c_{\rho,\mu}^\lambda$ as the number of connected components isomorphic to $\mg_\mu$ in $\mg_{\Skew}$. Applying Theorem~\ref{Skew}, any component that is isomorphic to $\mg_\mu$ has exactly one vertex $v_\mu$ whose row reading word is mapped to $U_\mu$ by $P$. Recall that SYam($\mu$) is defined to be the set of words $w$ such that $P(w)=U_\mu$.
Thus, $c_{\rho,\mu}^\lambda = |\{ w \in \textup{SYam}(\mu): w \textup{ is the row reading word of a filling in SYT}(\Skew)\}|$. \hfill $\square$

\vspace{.1in}
 
The next definition and corollary describe the effect of removing the signs and edges with the lowest labels from a signed colored graph.

\begin{dfn}\label{upward restriction}
\textup{Fix some signed colored graph $\mg=(V,\sigma, E_{m+1}\cup\dots\cup E_{n-1})$ and positive integer $h$. Let $\mh=(V,\sigma^\prime, E^\prime)$ be the signed colored graph defined by $\sigma^\prime_{i}=\sigma_{i+h}$ and $E^\prime_{i}=E_{i+h}$ whenever $\sigma_{i+h}$ and $E_{i+h}$ are defined in $\mg$, respectively. Then $\mh$ is termed \emph{the $h$-upward restriction of $\mg$}. If $h=1$, then $\mh$ is termed \emph{the upward restriction of $\mg$}. The restriction of $\mg$ to $(V, \sigma, E_{m+1}\cup\dots\cup E_{n-2})$ is termed \emph{the downward restriction of $\mg$}.
}
\end{dfn}

\begin{cor}\label{restrict up}
If $\mg$ is a dual equivalence graph, then the $h$-upward restriction of $\mg$ is a dual equivalence graph.
\end{cor}
\noindent $Proof:$ 
Without loss of generality, we may assume that $\mg=(V,\sigma, E_{m+1}\cup \ldots \cup E_{n-1})$ is connected. Let $\mh$ be the $h$-upward restriction of $\mg$. It follows immediately from Definition \ref{upward restriction} that $\mh$ obeys Axioms 1, 2, 3, 4, and 5, so we need only demonstrate Axiom 6. If $h < m$, then Definitions \ref{axioms} and \ref{upward restriction} guarantee that Axiom 6 holds for $\mg$ if and only if it holds for $\mh$. We may then restrict to the case where $\mg$ has type $(n,n)$ by simply taking the $(m-1)$-upward restriction. It thus suffices to consider $\mg=\mg_\lambda$ for some $\lambda \vdash n$, by Theorem~\ref{main}. 
Then $\mh$ can be obtained by removing the 1 through $h$ boxes from each vertex in $\mg_\lambda$ and subtracting $h$ from the values in all of the remaining boxes, creating some dual equivalence graph $\cup\mg_{\Skew}$, where the union is over all $\rho$ contained in $\lambda$ such that $\rho \vdash h$. By Theorem~\ref{Skew}, this union is a dual equivalence graph.
\hfill $\square$

\vspace{.1in} 

 The upward and the downward restrictions of a dual equivalence graph are structurally related, as is made precise in the following definition and corollary.

\begin{dfn}\label{color reversal}
\textup{Let $\mg=(V,\sigma,E)$, and let $\mh=(V,\sigma^\prime,E^\prime)$ be signed colored graphs such that
\begin{enumerate} \itemsep-.0in
\item $\sigma$ and $\sigma^\prime$ are maps onto $\{\pm 1\}^{N-1}$,
\item $\sigma_i=\sigma^\prime_{N-i}$ for all $1\leq i < N$,
\item $E_i=E^\prime_{N+1-i}$ whenever $E_i$ or $E^\prime_{N+1-i}$ is defined.
\end{enumerate}
Then $\mh$ is termed the \emph{color reversal} of $\mg$. 
}
\end{dfn}

\begin{cor}\label{inversion}
Let $\mg$ be an $(n,n)$-dual equivalence graph, and let $\mh$ be the color reversal of $\mg$. Then $\mg$ is isomorphic to $\mh$.
\end{cor}

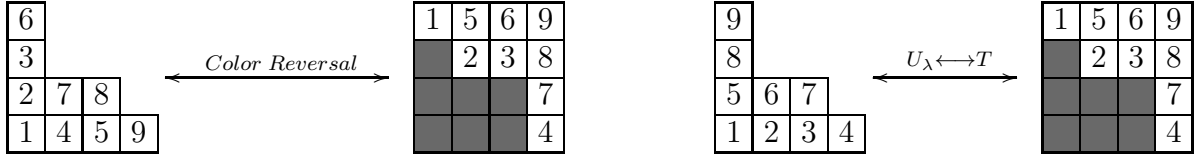
\begin{figure}[H]
   \begin{center} 
   	\begin{displaymath}
		\ytableausetup{aligntableaux=center}
		\xymatrix{
			\ytableaushort{6,3,278,1459} \ar@{<->}[rrr]^{Color \; Reversal} &&&  \;\;\ytableaushort{1569,\none 238,\none \none \none7,\none\none\none4}
			* [*(grey)]{0,1,3,3}
			}
	\hspace{.7in}
			\xymatrix{
			\ytableaushort{9,8,567,1234}  \ar@{<->}[rr]^{U_\lambda \longleftrightarrow T} &&  
			\;\; \ytableaushort{1569,\none 238,\none \none \none7,\none\none\none4}
			* [*(grey)]{0,1,3,3} }
\end{displaymath}
   \end{center}
 	 \caption{At left, a tableau and its color reversal. At right, $U_\lambda$ for $\lambda=(4,3,1,1)$ and its corresponding skew tableau $T\in \textup{SYT}(\Skewmu)$.}
 \label{color inversion}
\end{figure}

\noindent $Proof:$
It suffices to only consider connected graphs. Applying Theorem~\ref{main} allows us to further reduce to the case where $\mg$ is an arbitrary $\mg_\lambda$. Let $\mu/\rho$ be the skew shape given by rotating $\lambda$ by $180^\circ$, and let $\mg_{\Skewmu}$ be the dual equivalence graph on SYT($\mu/\rho$). The reader can check that $\mg_{\mu/\rho}$ is isomorphic to the color reversal of $\mg_\lambda$ by simply rotating any filling of $\lambda$ and then reversing the order of the numbers in the filling as in Figure \ref{color inversion}. In particular, $\mg_{\Skewmu}$ is connected. To show that $\mg_\lambda \cong \mg_{\mu/\rho}$, recall that $\mg_{\mu/\rho}$ is isomorphic to the signed colored graph induced by sending its vertices to their row reading words with edges given by $d_i$ and signature given by inverse descents. Applying Theorem~\ref{Skew}, $\mg_{\mu/\rho}\cong \mg_\lambda$ if the row reading word of any vertex $v$---and thus all vertices---in $\mg_{\mu/\rho}$ has sh$(v)=\lambda$.

Let $T\in \textup{SYT}(\Skewmu)$ be the skew tableau obtained by right justifying all of the rows of $U_\lambda$ and then top justifying all of the columns, as in the right side of Figure \ref{color inversion}. This transformation from $U_\lambda$ to $T$ is achieved by jeu de taquin, which preserves the shape of row reading words, as mentioned in Section \ref{The RSK Correspondence}. Therefore, the row reading word of $T$ has shape $\lambda$, completing the proof.
\hfill $\square$

\begin{remark} \textup{
Because color reversal acts as an isomorphism between dual equivalence graphs, it induces an isomorphism between standard dual equivalence graphs. While not proven here, it can be shown that this induced isomorphism acts on tableau via the Sch\"{u}tzenburger evacuation function as defined in \cite[Sec. 3.20]{Stanley}.}
\end{remark}

\subsection{Morphisms}

In this section, we set out to describe various properties of morphisms between signed colored graphs. We restrict our attention to graphs satisfying at least Axiom 1. For this reason, we may abuse notation and treat each $E_i$ as a function defined on those vertices admitting an $i$-edge. That is, we write $E_i(v)=w$ to mean that $v$ is contained in an $i$-edge with $w\neq v$. 

\begin{lem} \label{surjective}
Let $\mathcal{G}=(V,\sigma, E)$ and $\mathcal{H}=(V^\prime, \sigma^\prime, E^\prime)$ be nonempty signed colored graphs satisfying Axiom 1. Also, suppose that there exists a morphism $\phi\colon \mg \rightarrow \mh$. Then the following hold.
\begin{enumerate}\itemsep.0in
\item $\phi(E_i(v))=w$ if and only if $E^\prime_i(\phi(v))=w$.
\item If $\mathcal{H}$ is connected, then $\phi$ is surjective.
\item If $\phi$ is surjective and either $\mathcal{G}$ or $\mathcal{H}$ obey Axiom 2 or 3, they both do.
\item If $\phi$ is surjective and $\mathcal{G}$ obeys Axiom 4, 5, or 6, then $\mathcal{H}$ does as well.
\item If $\mh$ is connected, then $\phi$ is an $m$ to 1 map on vertices for some positive integer $m$.
\item If $\phi$ is a bijection from $V$ to $V^\prime$, then $\phi$ is an isomorphism from $\mg$ to $\mh$.
\end{enumerate}
\end{lem}

\noindent $Proof$: We begin with Part~1 and continue in order. Morphisms preserve signatures, so if either $v$ or $\phi(v)$ does not admit an $i$-edge, then Axiom 1 implies that neither is contained in an $i$-edge. Now suppose that $\{u,v\}$ and $\{\phi(v),w\}$ are $i$-edges in $\mg$ and $\mh$, respectively. By Axiom 1, these are the only $i$-edges containing $v$ and $\phi(v)$, so by the definition of morphisms, $\phi(u)=w$. That is, $\phi(E_i(v))= E^\prime_i(\phi(v))$. Thus, Part~1 holds.

For Part~2, we apply Part~1. Choose any $v \in \mg$. All of the vertices connected to $\phi(v)$ by an edge have a preimage in $\mg$ by Part~1. Since $\mh$ is connected, we may then induct to reach any vertex in $\mh$ by repeating this process.
 
Axioms 2 and 3 are concerned with the relationship between signatures and edges. By the definition of a morphism, $\phi$ preserves signatures, and by Part~1, $\phi(E_i(v))=w$ if and only if $E^\prime_i(\phi(v))=w$ for each $i$ where $E_i^\prime$ is defined. Hence, if either graph obeys Axiom 2 or 3, they both do.
  
 Proving Part~4 is a matter of repeatedly applying Part~1 to show that $\mh$ inherits the desired properties from $\mg$.
 As we will not need these properties for later results, we leave the details to the reader.

For Part~5, let $\phi^{-1}(x)$ be the set of vertices in $\mg$ mapped to $x$ by $\phi$. If $E_i^\prime(x)=y$, then Part~1 and Axiom 1 imply that the vertices in $\phi^{-1}(x)$ share an $i$-edge with distinct vertices in $\phi^{-1}(y)$ and vice versa. Hence, $|\phi^{-1}(x)| = |\phi^{-1}(y)|.$ Because $\mh$ is connected, we may induct from $x$ to any vertex in $\mh$ to show that the fiber over every vertex in $\mh$ has the same cardinality. That is, $\phi$ is an $m$ to one map on vertices.
   
Part~6 requires the existence of an inverse morphism, which follows from Part~1 and bijectivity.~\hfill~$\square$

\vspace{.1in}

\begin{remark}\label{covering}
\textup{
Let $\phi\colon \mg \rightarrow \mh$ be as in Lemma~\ref{surjective} and suppose that $\mh$ is connected. Then in the language of algebraic graph theory, Part~1 implies that $\phi$ is a covering map. Here we use the definition of a covering map on graphs - a surjective map that sends vertices to vertices and induces a bijection between edges containing $v$ and edges containing $\phi(v)$ for each vertex $v$ in $\mg$ - though the topological definition of a covering map can be made to apply as well. In this context, Parts 5 and 6 of Lemma~\ref{surjective} are well known properties of covering maps. For more details, see \cite[Sec. 6.8]{GR}.}
\end{remark}

\begin{cor}\label{terminal}
If $\phi$ is any morphism from a connected dual equivalence graph to a connected signed colored graph satisfying Axiom 1, then $\phi$ is an isomorphism.
\end{cor}
\noindent $Proof$: By Part~6 of Lemma~\ref{surjective}, we need only show that $\phi$ is bijective. By Parts 2 and 5 of Lemma~\ref{surjective}, $\phi$ is a surjective $m$ to one map, so it suffices to show that $m=1$. Taking restrictions if necessary, we may then assume that $\phi$ is a map between connected signed colored graphs of type $(n,n)$. Applying Theorem~\ref{main} allows us to assume that the domain of $\phi$ is some $\mg_\lambda$.  Now notice that in $\mg_\lambda$ there is only one vertex with the signature of $U_\lambda$, implying that $m=1$. Thus, $\phi$ is a bijection.  \hfill $\square$





\vspace{.1in}

In light of (\ref{schurs}) and 
Part~5 of Lemma~\ref{surjective}, it is natural to look for signed colored graphs that admit a morphism onto a union of standard dual equivalence graphs. The next theorem describes a class of signed colored graphs that admit such a morphism and gives a formula for their Schur expansion.

\begin{thm}\label{KnuthMorph}
Let $\mg=(V,\sigma, E)$ be an $(n,n)$-signed colored graph
satisfying the following properties:
\begin{enumerate}\itemsep-.0in
\item $\mg$ obeys Axiom 1,
\item The vertices in $V$ are indexed by a subset of $S_n$,
\item The signature function $\sigma$ is given by inverse descent sets of permutations,
\item $E_i(v)$ is Knuth equivalent to $d_i(v)$ for all $1<i < n$ and all $v \in V$ admitting an $i$-edge. 
\end{enumerate}

\noindent Then $P\colon V\rightarrow \textup{SYT}(n)$ induces a morphism $\phi\colon  \mg \rightarrow \bigcup_{\lambda \vdash n} \mg_\lambda$.
Furthermore, 
\[ \sum_{v\in V} F_{\sigma(v)} = \sum_{\lambda \vdash n} \, |\{V\cap \textup{SYam}(\lambda)\}| \cdot s_\lambda.
\]
\end{thm}

\noindent $Proof:$ To show that $P$ induces a morphism on $\mg$, first notice that Theorem~\ref{Skew} implies that $\sigma(P(w))=\sigma(w)$ for all $w \in V$, so $\phi$ preserves signatures. Now choose any $w\in V$ admitting an $i$-edge, and let $K_w$ be a composition of fundamental Knuth equivalences such that $E_i(w)=K_w\circ d_i(w)$. Because $P$ is constant on Knuth classes, $d_i(P(w))=d_i(P(K_w(w)))$. By Lemma~\ref{Knuth} and Theorem~\ref{Skew}, $d_i$ commutes with $K_w$ and $P$. Thus, 
\begin{equation*}
d_i(P(K_w(w)))= P(d_i\circ K_w(w))= P(K_w\circ d_i(w))= P(E_i(w)).
\end{equation*}
Hence, $\phi$ preserves edges and signatures, satisfying the definition of a morphism. 

To verify the second part of the theorem, we may restrict the domain to $\phi^{-1}(\mg_\lambda)$. Applying  (\ref{schurs}), we need only show that $|\{ V\cap \textup{SYam}(\lambda)\} |$ is equal to the value of the index $m$ of this restriction, as stated in Part~5 of Lemma~\ref{surjective}. By definition, SYam($\lambda$) is the set of words mapped to $U_\lambda$ by $P$, so $m=|\textup{SYam}(\lambda)\cap V|$, completing the proof.~\hfill~$\square$

\vspace{.1in}

Theorem~\ref{KnuthMorph} is close to a complete description of morphisms from $(n,n)$-signed colored graphs satisfying Axiom 1 onto dual equivalence graphs. This is made precise in the following theorem.

\begin{thm}\label{KnuthMorphCompleteness}
Let $\mg$ be an nonempty $(n,n)$-signed colored graph satisfying Axiom 1 and admitting a morphism $\phi\colon \mg \rightarrow \mg_\lambda$. Then $\mg$ is isomorphic to the $(n,n)$-restriction of some $\mh=(V^\prime, \sigma^\prime, E^\prime)$ satisfying the following properties:
\begin{enumerate}\itemsep.0in
\item $\mh$ obeys Axiom 1,
\item The vertices in $V^\prime$ are indexed by a subset of $S_N$ for some $N \geq n$,
\item The signature function $\sigma\colon V\rightarrow \{\pm1\}^{N-1}$ is given by inverse descent sets of permutations,
\item $E_i(v^\prime)$ is Knuth equivalent to $d_i(v^\prime)$ for all $1<i < n$ and all $v^\prime \in V^\prime$ admitting an $i$-edge. 
\end{enumerate}

\end{thm}

\noindent
\emph{Proof:} 
By Parts 2 and 5 of Lemma~\ref{surjective}, $\phi$ is a surjective $m$ to 1 map. Choose any $\mu \supset \lambda$ such that $|\textup{SYT}(\mu)| \geq m$. Here, $N=|\mu|$. Let $A$ be any skew tableau of shape $\mu/\lambda$ with the values $n+1,\ldots, N$. 
Let $\tilde \mg=(V, \sigma, E)$ be the $(n,N)$-signed colored graph with the same vertex and edge set as $\mg$ and with $\sigma(v)$ defined for all $v \in V$ as the signature of $\phi(v)$ augmented by $A$. It is clear that $\mg$ is the $(n,n)$-restriction of $\tilde \mg$, so we will construct $\mh$ to be isomorphic to $\tilde \mg$.

We will define $\mh$ by finding an appropriate relabeling of the vertices of $\tilde \mg$.  By the construction of $\tilde \mg$, we may extend $\phi$ to a morphism $\tilde \phi:\tilde \mg\rightarrow \mg_{\lambda,A}$. Because $|\textup{SYT}(\mu)|$ is greater than or equal to the index $m$ of $\tilde \phi$, there exists an injective map on the vertices in $V$, sending each vertex $v$ to $(\tilde \phi(v), T_{v})$ for some tableau $T_{v}$ of shape $\mu$. Since $\tilde \phi(v)$ and $T_{v}$ are the same shape, the inverse R-S-K correspondence takes each pair to a unique permutation in $S_N$. Let $f$ be the injective function taking vertices of $V$ to these permutations in $S_N$. We claim that if $\mh=(V^\prime,\sigma^\prime, E^\prime)$ is the signed colored graph induced by letting $f$ relabel the vertices of $\tilde \mg$, then $\mh$ satisfies all of the desired properties in the statement of the theorem. 

First, $\tilde \mg \cong \mh$, because $\tilde \mg$ and $\mh$ only differ by the labeling of their vertex sets. In particular, $\mh$ inherits Axiom 1 from $\tilde \mg$. Choose any $v \in V$ and $v^\prime \in V^\prime$ such that $f(v)=v^\prime$. To see that $\sigma^\prime$ agrees with the signature given by inverse descents, notice that $\tilde \phi(v) = P(v^\prime)$, by the construction of $f$. By the definition of a morphism, $\tilde \phi$ preserves signature, and by Theorem~\ref{Skew}, the signature of $P(v^\prime)$ is equal to the signature of $v^\prime$ given by inverse descents. Thus $\sigma^\prime$ is given by inverse descents. Lastly, we show that if $\{v, w\}$ is an $i$-edge in $\tilde \mg$, then $d_i(v^\prime)$ is in the Knuth class of $E^\prime_i(v^\prime)=f(w)$.  Applying Part~1 of Lemma~\ref{surjective} and Theorem~\ref{Skew} gives
\begin{equation*}
P(f(w))=\tilde \phi(w)= d_i(\tilde \phi(v^\prime)) = d_i(P(v^\prime))=P(d_i(v^\prime)),
\end{equation*}
and so $f(w)$ is in the Knuth class of $d_i(v^\prime)$.~\hfill~$\square$

\vspace{.1in}

\begin{remark}\textup{
Theorems \ref{KnuthMorph} and \ref{KnuthMorphCompleteness} can both be extended to statements about $(n,N)$-signed colored graphs with morphisms to augmented dual equivalence graphs.}
\end{remark}

\subsection{Local Conditions for Axiom 6}\label{ax6 equiv}

Out of the six dual equivalence axioms, Axiom 6 is the only one that cannot be checked by testing local criteria. In this section we will show that an equivalent axiomatization is given by strengthening Axiom 4 and omitting Axiom 6.

\begin{dfn}\label{def 4plus} \textup{ 
A signed colored graph $\mg=(V, \sigma, E_{m+1}\cup\ldots\cup E_{n-1})$ is said to obey Axiom $4^+$ if for all $m+1 < i < n$, any component of the edge colored graph $(V, E_{i-3}\cup E_{i-2} \cup E_{i-1} \cup E_{i})$ is isomorphic to a component of the restriction of some $\mg_\lambda=(V^\prime, \sigma^\prime, E^\prime)$ to $(V^\prime,E^\prime_{i-3} \cup E^\prime_{i-2}\cup E^\prime_{i-1} \cup E^\prime_{i})$, where $E_{i-3}$ or $E_{i-2}$  is omitted if $i \leq m+3$ or $i\leq m+2$, respectively.
}
\end{dfn}

\noindent We now state the main result of this section. The proof is postponed until after a necessary lemma.

\begin{thm}\label{4plus}
A signed colored graph satisfies Axioms 1, 2, 3, $4^+,$ and 5 if and only if it is a dual equivalence graph.
\end{thm}

\begin{remark}\label{relabeling 4}\textup{
\vspace{-.1in}
\begin{enumerate}\itemsep-.2in
\item
We may readily classify the set of edge colored graphs described in Definition \ref{def 4plus}, i.e., the set of edge colored graphs that arise as  components of the restriction of some $\mg_\lambda=(V, \sigma, E)$ to $(V, E_{i-3}\cup E_{i-2}\cup E_{i-1} \cup E_{i})$ for all choices $1<i<|\lambda|$.  Applying Corollary~\ref{restrict up}, each such edge colored graph is the result of choosing an appropriate $\lambda\vdash 6$, restricting $\mg_\lambda=(V^\prime, \sigma^\prime, E^\prime)$ to the edge colored graph $(V^\prime, E^\prime)$, and adding a fixed nonnegative integer to each edge label. \\
\item
\indent A signed colored graph $\mg =  (V, \sigma, E_{m+1} \cup \ldots \cup E_{n-1})$ satisfies Axioms 1, 2, 3, and $4^+$ if and only if, for any $m< i < n$, each component of $(V, \sigma,  E_{i-3}\cup E_{i-2} \cup E_{i-1} \cup E_{i})$ is isomorphic to a component of the restriction of some $\mg_\lambda=(V^\prime, \sigma^\prime, E^\prime)$ to $(V^\prime, \sigma^\prime, E_{i-3}\cup E_{i-2} \cup E_{i-1} \cup E_{i})$, where $E_{i-3}, E_{i-2}$, or $E_{i-1}$ is omitted if $i\leq m+3, i\leq m+2$ or $i=m+1$, respectively.
\end{enumerate}
}
\end{remark}

\begin{lem}\label{connectedness condition}
Let $\mg$ be a connected $(n+1,n+1)$-signed colored graph satisfying Axioms 1, 2, 3, 4, and 5 whose downward restriction and upward restriction are both dual equivalence graphs. Let $\mc$ be any component of the downward restriction of $\mg$. Additionally, suppose that for every pair of distinct $(n,n+1)$-components $\ma$ and $\mb$ that are neighbors of $\mc$ in $\mg$, there exists a component of the upward restriction of $\mg$ whose vertex set intersects $\ma$ nontrivially and intersects $\mb$ nontrivially. Then $\mg$ is a dual equivalence graph. 
\end{lem}
\noindent
\emph{Proof:} 
Consider any such $\ma$ and $\mb$ whose vertices are nontrivially intersected by some component $\md$ of the upward restriction of $\mg$. Because $\md$ obeys Axiom 6, there exists a path in $\mg$ from a vertex in $\ma$ to a vertex in $\mb$ crossing a single $n$-edge, which would be labeled as an $(n-1)$-edge in $\md$. Hence, $\ma$ and $\mb$ are neighbors in $\mg$. Because $\ma$ and $\mb$ were chosen arbitrarily, all neighbors of $\mc$ are pairwise neighbors of each other in $\mg$.

Next we show that every pair of vertices in $\mg$ can be connected by a path containing at most one $n$-edge. By Corollary~\ref{complete}, every $(n, n+1)$-component of $\mg$ has the same number of neighbors in $\mg$. Thus, if $\mb$ is any neighbor of $\mc$ in $\mg$, each neighbor of $\mb$ is either $\mc$ or a neighbor of $\mc$ in $\mg$. That is, 
$\mc \cup (\bigcup \mb)$ contains all of the vertices of $\mg$, where the union is over all $\mb$ that neighbor $\mc$ in $\mg$.
Thus, all $(n, n+1)$-components of $\mg$ are pairwise neighbors. In particular, any two vertices in $\mg$ can be connected by a path crossing at most one $n$-edge.

This property of paths in $\mg$, along with the hypothesis that the $(n, n+1)$-restriction of $\mg$ is a dual equivalence graph, guarantees that $\mg$ satisfies Axiom 6. By assumption, $\mg$ also satisfies Axioms 1, 2, 3, 4, and 5, so $\mg$ is a dual equivalence graph.~\hfill~$\square$

\vspace{.2in}

\noindent
\emph{{\bf Proof of Theorem~\ref{4plus}:}}

  We begin with the forward implication by showing that any signed colored graph satisfying Axioms $1, 2, 3, 4^+$, and 5 must also satisfy Axiom 6.
Because Axiom 6 is concerned with edge sets, it suffices to only consider signed colored graphs of type $(n,N)$, as in the proof of Corollary~\ref{restrict up}. We now proceed by induction on $n$. Axiom $4^+$, when considered with Axioms 1, 2, and 3, implies the theorem holds for $(n,N)$-signed colored graphs with $n \leq 6$. Now suppose that the result holds for all $(n,N)$-signed colored graphs, and consider any $(n+1, n+1)$-signed colored graph obeying axioms 1, 2, 3, $4^+$, and 5.  Call this graph $\mg$.

Choose any component $\mc$ of the downward restriction of $\mg$. We will show that $\mg$ and $\mc$ satisfy the hypotheses of Lemma~\ref{connectedness condition}. By assumption, $\mg$ satisfies Axioms 1-5, and by induction, the upward and downward restrictions of $\mg$ are dual equivalence graphs. It then remains to be shown that if $\ma$ and $\mb$ are any distinct neighbors of $\mc$ in $\mg$, then there exists a component $\md$ of the upward restriction of $\mg$ that intersects the vertices of $\ma$ and $\mb$ nontrivially.

Next, we  label the vertices in $\mg$ by tableaux. By Theorem~\ref{main}, $\mc$ is isomorphic to some $\mg_\lambda$, and by Proposition \ref{extension} there is some $\mg_\mu$ isomorphic to the unique extension of $\mc$. Furthermore, Corollary~\ref{complete} guarantees the existence of a morphism from $\mg$ to $\mg_\mu$. Label the vertices of $\mg$, as well as the vertices of the downward restrictions of $\mg$, by the set of tableaux in SYT($\mu$) as given by the image of this morphism. By Corollary~\ref{terminal},  this morphism restricts to an isomorphism from any given component of the downward restriction of $\mg$ to some component of the downward restriction of $\mg_\mu$.

\begin{figure}[H]
\begin{center}
\scalebox{.85}{$\begin{array}{cl}
\xymatrix{
	\ytableaushort{34,125}\; \ar@{=}[r]^{2}_{3} & \;\ytableaushort{24,135} }
&\!\! \xymatrix{
	 \ar@{->}[rr] &&\;\;\;\; \ytableaushort{{}{},{}{}5} }\\ & \\
\xymatrix{
	\ytableaushort{25,134}\; \ar@{-}[r]^{2} & \;\ytableaushort{35,124}\; \ar@{-}[r]^{3} & \;\ytableaushort{45,123}} 
&\!\! \xymatrix{
 \ar@{->}[rr] &&  \;\;\;\;\ytableaushort{{}5,{}{}{}} 
			 }
\end{array}$}
\end{center}
\caption{The $(4,4)$-components of $\mg_{(3,2)}$, at left, are represented via the cells containing 5, at right. Signatures are omitted.}\label{Type by Corner}
\end{figure}
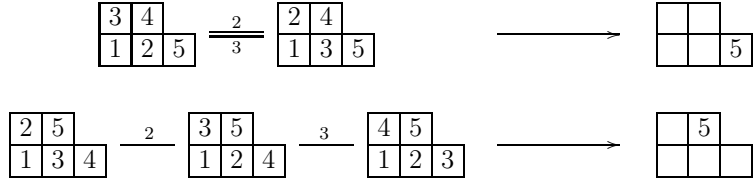

We may now associate the isomorphism types of $\ma, \mb,$ and $\mc$ to the position of the cell containing $n+1$ in fillings of $\mu$, as in Figure~\ref{Type by Corner}. It follows from Theorem~\ref{Haiman} that the number of $(n,n+1)$-components in $\mg_\mu$ is equal to the number of Northeast corners of $\mu$, with the isomorphism type of each $(n, n+1)$-component determined by which Northeast corner is filled by $n+1$. By Corollary~\ref{complete}, $\ma$ can be described as the unique neighbor of $\mc$ whose tableaux have $n+1$ in a particular Northeast corner of $\mu$. Since $\ma, \mb$, and $\mc$ have distinct isomorphism types, they each have different Northeast corners $c_1, c_2$, and $c_3$, respectively, filled by $n+1$.

To apply Lemma~\ref{connectedness condition}, it suffices to show that there is a tableau $T \in \textup{SYT}(\mu)$ from which $n+1$ may be moved into each of these three Northeast corners while staying in the vertex set of a component of the upper restriction,  i.e., without crossing a 2-edge.
It is possible to use jeu de taquin to describe which cells $n+1$ may be moved to without crossing a 2-edge. Let $S \in \textup{SYT}(\mu)$ be some vertex of $\mg_\mu$. As in the proof of Corollary~\ref{restrict up}, the set of vertices in the same component of the upward restriction of $\mg_\mu$ as $S$ is described by the component of $S$ after omitting the 1-cell (for notational simplicity, we choose not to relabel the boxes). We may then perform jeu de taquin to retrieve some unique tableau $J(S)$ such that sh$(J(S))$ is a straight shape. 
Any sequence of jeu de taquin slides acts on the row reading words of a filling by a sequence of fundamental Knuth equivalences, and so commutes with the action of $d_i$ by Lemma~\ref{Knuth}. Theorem~\ref{Haiman} guarantees that $n+1$ may then be moved to any Northeast corner of sh$(J(S))$. In particular, if $c$ is some Northeast corner of $\mu$, $x$ is the value assigned to $c$ by $S$, and $x$ is also in a Northeast corner of $J(S)$, then there is a path in $\mg_\mu$ with no 2-edges that connects $S$ to a tableau with the value $n+1$ in $c$. Therefore, to apply Lemma~\ref{connectedness condition}, we need only find $T \in \textup{SYT}(\mu)$ such that the values in $c_1$, $c_2$, and $c_3$ are in Northeast corners of $J(T)$. 

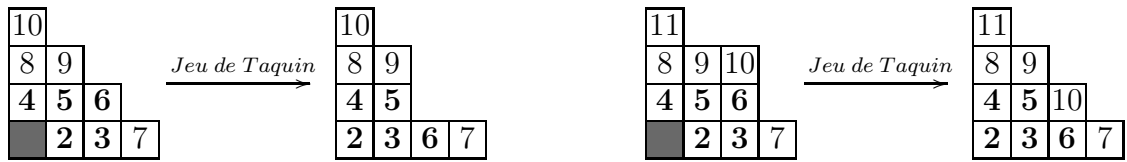
\begin{figure}[h]
   \begin{center} 
   	\begin{displaymath}
		\xymatrix{
			\ytableaushort{{10},89,{\bf{4}}{\bf{5}}{\bf{6}},{*(grey){}}{\bf{2}}{\bf{3}}7} \ar@{->}[rr]^{Jeu \;de\; Taquin } && \;\; \ytableaushort{{10},89,{\bf{4}}{\bf{5}},{\bf{2}}{\bf{3}}{\bf{6}}7}  && 
			\ytableaushort{{11},89{10},{\bf{4}}{\bf{5}}{\bf{6}},{*(grey){}}{\bf{2}}{\bf{3}}7} \ar@{->}[rr]^{Jeu\; de \; Taquin } && \;\;\ytableaushort{{11},89,{\bf{4}}{\bf{5}}{10},{\bf{2}}{\bf{3}}{\bf{6}}7}
			}
\end{displaymath}
   \end{center}
 	 \caption{At left, a filling chosen to lose the corner containing 6 to jeu de taquin. At right, a filling of a nonstaircase shape that does not lose any Northeast corners to jeu de taquin.}
 \label{JDT}
\end{figure}

Having reduced the problem to a matter of jeu de taquin, we now consider two cases. First, suppose $\mu$ has at least four Northeast corners. Given $c_1$, $c_2$, and $c_3$, we are free to lose some fourth corner $c_4$.  By filling all boxes weakly southwest of $c_4$ with as low of values as possible, the fourth corner will always be moved, leaving the three given corners unchanged, as in the left half of Figure \ref{JDT}.

As a second case, suppose $\mu$ is not a staircase shape, i.e. $\mu$ is not of shape $(\mu_1, \mu_1 - 1, \ldots, 2, 1)$. Then we claim that there exists some $T\in \textup{SYT}(\mu)$ such that every value in a Northeast corner of $T$ is in a Northeast corner of $J(T)$. Notice that $\mu$ must contain some rectangle whose Northeast corner is not a Northeast corner of $\mu$ but is either the most east cell in its row or the most north cell in its column. It is easy to check that any filling of $\mu$ such that all values inside of the rectangle are less than all values outside of the rectangle will suffice. See the right half of Figure \ref{JDT} for an example.

We have shown that $\mg$ is a dual equivalence graph if $\mu$ is not a staircases or has at least four Northeast corners. There are only three staircase shapes that have strictly less than four Northeast corners, all of which have size less than or equal to 6. All shapes with size less than or equal to 6 are contained in the base case, so $\mg$ is a dual equivalence graph.

Assuming our inductive hypotheses, we have demonstrated that any $(n+1,n+1)$-signed colored graph that satisfies Axioms 1, 2, 3, $4^+$, and 5 also satisfies Axiom 6. An $(n+1,N)$-signed colored graph obeys Axiom 6 if and only if its $(n+1, n+1)$-restriction obeys Axiom 6, completing the inductive step. Hence, any signed colored graph obeying Axioms 1, 2, 3, $4^+$, and 5 is a dual equivalence graph.

The reverse implication follows more quickly. If $\mg=(V,\sigma, E)$ is a dual equivalence graph, then we need only show that $\mg$ obeys Axiom $4^+$. As in Part~1 of Remark~\ref{relabeling 4}, we may assume that $\mg$ has type $(n,n)$. By Theorem~\ref{main}, $\mg$ is isomorphic to some $\mg_\lambda$, so axiom $4^+$ follows immediately. Hence, every dual equivalence graph satisfies Axioms 1, 2, 3, $4^+$, and 5, completing the proof.~\hfill~$\square$

\vspace{.1in}

\begin{figure}[H]
   \begin{center} 
   \begin{displaymath}
	\scalebox{.85}{
	\xymatrix{
		&\star \ar@{-}[r]^4&\star\ar@{-}[r]^2&\star\ar@{=}[r]^3_4&\star\ar@{-}[r]^5&\bullet\ar@{-}[r]^3&\bullet&&\\
		\bullet\ar@{-}[ur]^5\ar@{=}[dr]^2_3&&\star\ar@{=}[ul]^3_2\ar@{-}[dl]_5&&&\bullet\ar@{-}[ur]^2\ar@{=}[dr]^4_5&&\bullet\ar@{=}[ul]^5_4\ar@{-}[dl]_2&\\
		&\bullet\ar@{-}[r]^4&\bullet\ar@{-}[r]^2&\bullet\ar@{=}[r]^3_4&\bullet&\bullet\ar@{-}[r]^3&\bullet&&\\
		&&\bullet\ar@{-}[r]^3&\bullet\ar@{-}[ur]^<<<<<5&\bullet\ar@{=}[r]^3_4\ar@{-}[ur]^<<<<<5&\bullet\ar@{-}[r]^2&\bullet\ar@{-}[r]^4&\bullet&\\
		&\bullet\ar@{-}[ur]^2\ar@{=}[dr]^4_5&&\bullet\ar@{=}[ul]^5_4\ar@{-}[dl]_2&&&\bullet\ar@{=}[dr]^2_3\ar@{-}[ur]^5&&\bullet\ar@{=}[ul]^3_2\ar@{-}[dl]_5\\
		&&\bullet\ar@{-}[r]^3&\bullet\ar@{-}[r]^5&\bullet\ar@{=}[r]^3_4&\bullet\ar@{-}[r]^2&\bullet\ar@{-}[r]^4&\bullet&
			}}
\hspace{.2in}
		\scalebox{.85}{
		\xymatrix{&&\\&&\\
			{\young(6,\hfill\hfill,\hfill\hfill\hfill)} \;\; \ar@{-}[r]^{5} \ar@{-}[d]^<<<<<{5} & \;\; {\young(\hfill,\hfill\hfill,\star\hfill6)}\ar@{-}[r]^{5} &{\young(\hfill,\hfill 6,\hfill\hfill\hfill)}   \ar@{-}[d]^>>>>>{5}\\
			 {\young(\hfill,\hfill 6,\hfill\hfill\hfill)} \ar@{-}[r]^{5} & \;\; {\young(\hfill,\hfill\hfill,\hfill\hfill 6)} \ar@{-}[r]^{5}&{\young(6,\hfill\hfill,\hfill\hfill\hfill)}}
			}
	\end{displaymath}   
   \end{center}
 	 \caption{Two representations of the smallest graph in $\mathcal{F}$ with edge labels in $\{2, 3, 4, 5\}$. Starred vertices on the left correspond to the starred shape on the right.}
 \label{smallest example}
\end{figure}
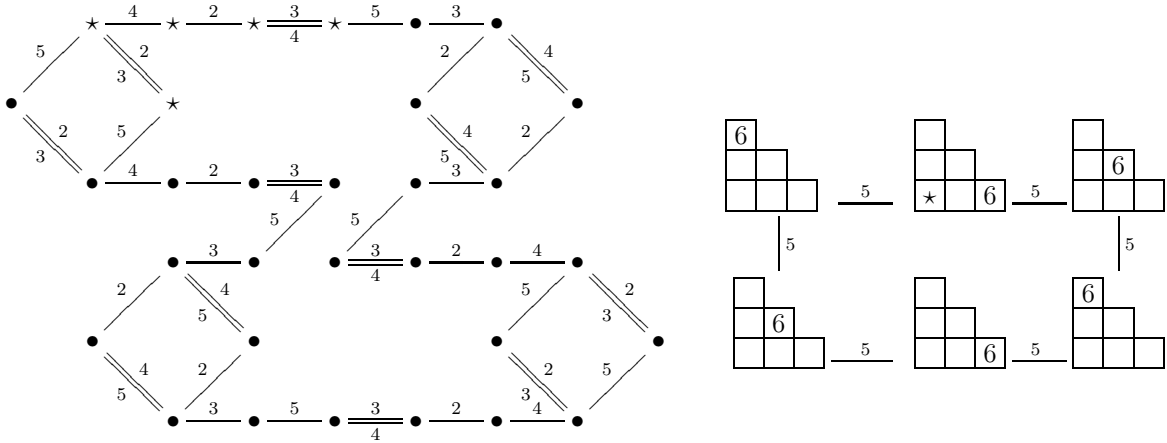

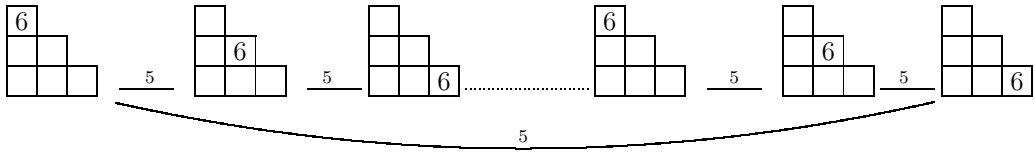
\begin{figure}[h]
   \begin{center} 
   \begin{displaymath}
	\scalebox{.85}{
	\xymatrix{
		{\young(6,\hfill\hfill,\hfill\hfill\hfill)}\;\; \ar@{-}[r]^{5} \ar@{-}@/_2.2pc/[rrrrrr]^{5} & \;\; {\young(\hfill,\hfill6,\hfill\hfill\hfill)}\;\; \ar@{-}[r]^{5}&{\young(\hfill,\hfill\hfill,\hfill\hfill6)} \ar@{.}[rr]&&{\young(6,\hfill\hfill,\hfill\hfill\hfill)}\;\; \ar@{-}[r]^{5}& \;\; {\young(\hfill,\hfill6,\hfill\hfill\hfill)}\ar@{-}[r]^{5} &{\young(\hfill,\hfill\hfill,\hfill\hfill6)} 
		}}
	\end{displaymath}
   \end{center}
   \caption{A generic element of $\mathcal{F}$ in the case where edge labels are in $\{ 2, 3, 4, 5\}$.}
 \label{F}
\end{figure}

We have actually proven a slightly stronger---and sometimes easier to check---condition. In the previous proof, Axiom $4^+$ was only invoked when considering the staircase with six cells, while the usual Axiom 4 could have been used for the staircase with three cells. If we assume Axioms 1, 2, 3, 4, and 5, then the staircase with six cells can only break Axiom 6 in a specific set of graphs. The (5,6)-components of $\mg_{(3,2,1)}$ have three distinct isomorphism types. By Corollary~\ref{complete}, each of these components is connected to two neighboring components. If Axiom 6 is not satisfied, then Lemma~\ref{connectedness condition} implies these two components cannot be neighbors of each other. Rather, the (5,6)-components must form a loop. The smallest example is presented in Figure \ref{smallest example}. In the figure, there are two copies of each (5,6)-component. For each positive integer $m\geq 2$, there is then a unique graph with $m$ isomorphic copies of each (5,6)-components, as in Figure \ref{F}. We may then omit signatures and relabel edges with elements in $\{i-3, i-2, i-1, i\}$, as mentioned in Remark~\ref{relabeling 4} for a full description of how Axiom $4^+$ can break in the presence of Axioms 1, 2, 3, 4, and 5. 

Let $\mathcal{F}$ be the set of edge colored graphs with edge sets $E_{i-3}\cup\ldots\cup E_{i}$ that satisfy Axioms 1, 2, 3, 4, and 5 but not Axiom 6. The following corollary reformulates Theorem~\ref{4plus} in terms of $\mathcal{F}$.

\begin{cor}\label{computer check}
Let $\mg=(V,\sigma, E_{m+1}\cup\ldots\cup E_{n-1})$ be a signed colored graph satisfying Axioms 1, 2, 3, 4, and 5. Then $\mg$ is a dual equivalence graph if and only if for all $m+4 < i < n$, the restriction of $\mg$ to the edge colored graph $(V, E_{i-3}\cup E_{i-2} \cup E_{i-1} \cup E_{i})$ has no components isomorphic to an element of $\mathcal{F}$.
\end{cor}

\begin{remark}
\textup{For any edge colored graph in $\mathcal{F}$, every vertex shares an edge with at least two other vertices. We may then give yet another characterization of dual equivalence graphs.  Let $\mg=(V, \sigma, E_{m+1}\cup \ldots \cup E_{n-1})$ be a signed colored graph obeying Axioms 1, 2, 3, 4, and 5. Choose $\mc$ to be any component of the restriction of $\mg$ to the edge colored graph $(V, E_{i-3} \cup E_{i-2} \cup E_{i-1} \cup E_{i})$ such that $m+3< i <n$ and the vertices of $\mc$ all have at least two adjacent vertices in $\mc$. Then $\mg$ is a dual equivalence graph if and only if $\mc$ is not in $\mathcal{F}$ for any choice of $\mc$. This characterization of dual equivalence graphs is used in the computer verification of Theorem~\ref{LLT graph}.
}
\end{remark}

\section{LLT and Macdonald Polynomials}\label{LLTs}

In this section we show that a family of LLT polynomials can be generated as a sum over vertices of a dual equivalence graph, and we provide a simple Schur expansion for polynomials in this family. Along the way, we recall Assaf's set of signed colored graphs associated to LLT polynomials and classify which ones are dual equivalence graphs.

\subsection{The Schur Expansion of $\llt$ when diam($\bnu) \leq 3$}

Recall the notation and vocabulary introduced in Section \ref{Symmetric Functions} on LLT polynomials and Macdonald polynomials.
To state the main theorem of this section, we will also need the following definition.


\begin{dfn}\label{diameter}
\textup{Given a $k$-tuple of skew shapes $\bnu$, let $S(\bnu)$ be the set of distinct shifted contents of the cells in $\bnu$. Define the} diameter \textup{of $\bnu$, denoted diam($\bnu$), as} 
\[
\textup{diam($\bnu$) := max$\{|R|: R\subset S(\bnu)$ and $|x-y|\leq k$ for all $x,y\in R\}$.} 
\]
\end{dfn}

\begin{figure}[H]
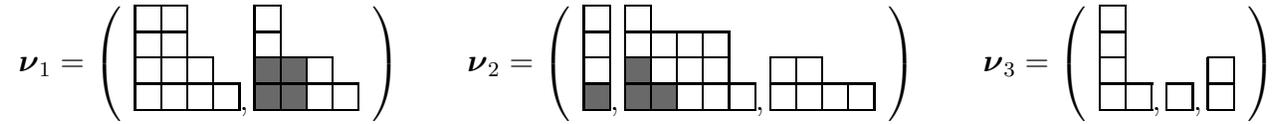

   \[  \ytableausetup{aligntableaux=bottom, boxsize=.33cm}
   \bnu_1=\left(
   \begin{array}{c} \ydiagram{2,2,3,4}, \ydiagram{1,1,2+1,2+2}*[*(grey)]{0,0,2,2} 
   \end{array} \right) \quad \quad
   \bnu_2=\left(
   \begin{array}{c}\ydiagram{1,1,1,0}*[*(grey)]{0,0,0,1}, \ydiagram{1,4,1+3,2+3}*[*(grey)]{0,0,1,2}, \ydiagram{0,0,2,4,}
   \end{array} \right) \quad\quad
   \bnu_3=\left(
   \begin{array}{c}\ydiagram{1,1,1,2}, \ydiagram{0,0,0,1}, \ydiagram{0,0,1,1}
   \end{array} \right)
   \]
  \caption{The tuples $\bnu_1$ and $\bnu_2$ have diameter 3. The tuple $\bnu_3$ has diameter 4.}
 \label{diam 3 tuples}
\end{figure}

\noindent See Figure \ref{diam 3 tuples} for an example.

\begin{remark}\label{k plus one}
\textup{
If $\bnu=(\nu^{(0)}, \ldots, \nu^{(k-1)} )$ is a tuple of skew shapes, it follows from Definition \ref{diameter} that diam$(\bnu)\leq k+1$. There are many examples where this bound is sharp. For instance, if $\bnu$ is a $k$-tuple of straight shapes where  $\nu^{(0)}$ has as least two columns and $|\nu^{(i)}|\geq 1$ for all $1 \leq i \leq k-1$, then diam$(\bnu)=k+1$. In Figure~\ref{diam 3 tuples}, $\bnu_3$ is an example of such a tuple.
}
\end{remark}

To further ease the presentation of the following results, we define two sets,
\begin{eqnarray}
&{\bf T}(\bnu, \lambda):=\{ {\bf T} \in \textup{SYT}(\bnu): \textup{the shifted content word of } {\bf T} \textup{ is in } \textup{SYam}(\lambda)\},\\
&T(\Skewmu, \lambda):=\{ T\in\textup{SYT}(\Skewmu): \textup{the content reading word of } T \textup{ is in } \textup{SYam}(\lambda)\}.
\end{eqnarray}
We are now able to present the main theorem of this section.

\begin{thm}\label{LLT Decomposition}
Let $\bnu$ be any $k$-tuple of skew shapes with \textup{diam}$(\bnu) \leq 3$. Then
\[
\llt = \sum_{\lambda \vdash |\bnu|} \;
\sum_{{\bf T} \in {\bf T}(\bnu, \lambda)} \hspace{-.1in} q^{\textup{inv}_k({\bf T})}s_\lambda .
\]
\end{thm}

\noindent The proof of Theorem~\ref{LLT Decomposition} is postponed until Section \ref{proof of LLT}. As mentioned in Remark~\ref{k plus one}, the set of $\bnu$ such that diam$(\bnu) \leq 3$ properly contains the set of $\bnu$ that are 2-tuples.  

The next corollary follows immediately by applying Theorem~\ref{LLT Decomposition} to the definition of modified Macdonald polynomials in (\ref{Macdonalds}).  We also use the fact that tuples of ribbons in TR($\Skewmu$) have diameter less than or equal to
three if and only if $\Skewmu$ does not contain (3,3) or (4) as a subdiagram. This fact is easily shown by noticing that the cells of every tuple of ribbons in TR($\Skewmu$) have the same set of distinct shifted contents.

\begin{cor}\label{Mac Decomposition}
Let $\mu/\rho$ be a skew shape not containing $(3,3)$ or $(4)$ as a subdiagram. Then
\[
\Mac = \sum_{\lambda \vdash |\Skewmu|} \;
\sum_{T \in T(\Skewmu, \lambda)}  q^{\textup{inv}(T)}t^{\textup{maj}(T)}s_\lambda.
\]
\end{cor}

\noindent In particular, Corollary~\ref{Mac Decomposition} applies to all $\widetilde H_{\mu}(X;q,t)$ with $\mu_1\leq3$ where $\mu_2 \leq 2$. 

\begin{remark}
\textup{
The conditions on $\bnu$ and $\mu/\rho$ in Theorem~\ref{LLT Decomposition} and Corollary~\ref{Mac Decomposition}, respectively, are sharp in the following sense. 
\textup{Let $\lambda = (2,2)$ and $\bnu = ((2),(1),(1))$ or 
$((1), (1), (1), (1))$. In particular, diam($\bnu)=4$, and $\bnu$ has diagram
\begin{center}
	\scalebox{.85}
	{$\left( \begin{array}{c} {\yng(2)}, { \yng(1)}, { \yng(1)} \end{array} \right)$}
 	\quad \textup{ or } \quad
	\scalebox{.85}
	{$\left( \begin{array}{c} {\yng(1)}, { \yng(1)}, { \yng(1)}, { \yng(1)} \end{array} \right)$}.
\end{center}
Then}
}
\begin{equation}
\llt \big|_{q^2s_{\lambda}}=  1  \quad \textup{and } \quad
 \sum_{{\bf T} \in {\bf T}(\bnu, \lambda)}  q^{\textup{inv}_k(T)}s_\lambda \big|_{q^2s_{\lambda}}= 0. \hspace{.15in}
\end{equation}
\textup{If $\lambda=(2,2)$ and $\mu = (4)$, then}
\begin{equation}
\widetilde H_{\mu}(X;q,t) \big|_{q^2s_{\rho}}=  1 \quad \textup{ and } \quad
 \sum_{T \in T(\mu, \lambda)}  q^{\textup{inv}(T)}t^{\textup{maj}(T)}s_\lambda \big|_{q^2s_{\lambda}}= 0.
\end{equation}
\textup{If $\lambda = (4,2)$ and $\mu = (3,3)$, then}
\begin{equation}
\widetilde H_{\mu}(X;q,t)\big|_{q^2s_{\rho}}=  1 \quad \textup{ and } \quad
 \sum_{T \in T(\mu, \lambda)}  q^{\textup{inv}(T)}t^{\textup{maj}(T)}s_\lambda \big|_{q^2s_{\lambda}}= 0.
\end{equation}

\end{remark}

\begin{cor}\label{Mac Decomposition2}
Let $\mu/\rho$ be a skew shape not containing $(2,2,2)$ or $(1,1,1,1)$ as a subdiagram. Then
\[
\Mac = \sum_{\lambda \vdash |\Skewmu|} \;
\sum_{T \in T(\tilde\mu/\tilde\rho, \lambda)}  q^{\textup{maj}(T)}t^{\textup{inv}(T)}s_\lambda.
\]
\end{cor}

\noindent \emph{Proof:} The corollary follows from (\ref{Mac conjugate}) and Corollary~\ref{Mac Decomposition}.~\hfill~$\square$
\vspace{.1in}

\noindent In particular, Corollary~\ref{Mac Decomposition2} applies to all $\tilde H_{\mu}(X;q,t)$ where $\mu$ has at most three rows and $\mu_3\leq1$.

\subsection{LLT graphs}\label{proof of LLT}

The goal of this section is to prove Theorem~\ref{LLT Decomposition}. We begin by following \citep{Assaf} in defining an involution that will provide the edge sets of a signed colored graph. In this section, $\bnu$ will always denote a $k$-tuple of skew shapes whose sizes sum to $|\bnu|=n$. Also, $w$ will always denote a permutation in $S_n$. 

Let 
 the involution $\tilde d_i:S_n \rightarrow S_n$ act by permuting the entries $i-1, i,$ and $i+1$ as defined by,
\begin{equation}
\begin{split}
\tilde d_i(\ldots i-1 \ldots i \ldots i+1 \ldots) = (\ldots i-1 \ldots i \ldots i+1 \ldots), \\
\tilde d_i(\ldots i+1 \ldots i \ldots i-1 \ldots) = (\ldots i+1 \ldots i \ldots i-1 \ldots), \\
\tilde d_i(\ldots i \ldots i-1 \ldots i+1 \ldots) = (\ldots i-1 \ldots i+1 \ldots i \ldots), \\
\tilde d_i(\ldots i \ldots i+1 \ldots i-1 \ldots) = (\ldots i+1 \ldots i-1 \ldots i \ldots).
\end{split}
\end{equation}
\noindent For instance, $\tilde d_3 \circ \tilde d_2(4123) = \tilde d_3(4123) = 3142.$ 

\begin{figure}[H]
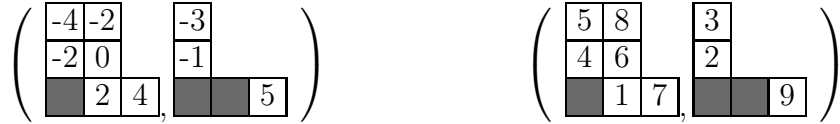

  \ytableausetup{boxsize=.19in}
  \[ \left( \begin{array}{c} \ytableaushort{{-4}{-2},{-2}0,{*(grey){}}24}, \ytableaushort{{-3},{-1},{*(grey){}}{*(grey){}}5}
   \end{array} \right) \hspace{1in}
   \left( \begin{array}{c}  \ytableaushort{ 58,46,{*(grey){}}17}, \ytableaushort{ 3,2,{*(grey){}}{*(grey){}}9}
   \end{array} \right)
  \]
  \caption{On the left, the shifted contents of a pair of skew diagrams with $\tau = 456667899$. On the right, a standard filling of the same tuple with shifted content word 534826179.}
\label{tau}
\end{figure}

To decide when to apply $d_i$ and when to use $\tilde d_i$, we appeal to the shifted content. Numbering the cells of a fixed $\bnu$ from 1 to $n$ in shifted content reading order, let $\tilde c_i$ be the shifted content of the $i^{th}$ cell. Define the weakly increasing word $\tau =\tau_1 \tau_2 \ldots \tau_n$ by 
\begin{equation}
\tau_i = \textup{max}\{j\in [n]: \tilde c_j - \tilde c_i \leq k\}.
\end{equation}
\noindent See Figure \ref{tau} for an example. To emphasize the relationship between $\tau$ and $\bnu$, we will sometimes write $\tau=\tau(\bnu)$. Notice that there are finitely many possible $\tau$ of any fixed length $n$. Specifically, $\tau$ will always satisfy $\tau_n = n$ and $i\leq \tau_i \leq \tau_{i+1}$ for all $i < n$. In fact, the number of possible $\tau$ is the $n^{th}$ Catalan number (for details on the Catalan numbers and an extensive list on where they arise in mathematics, see \citep{Catalan}). 
Next, let 
  $m(i)$ be the index of the the value in $\{i-1, i, i+1\}$ that occurs first in $w$, and let  $M(i)$ be the index of the value in $\{i-1, i, i+1\}$ that occurs last  in $w$. We now define the desired involution,
\begin{equation} 
D_i^{(\tau)}(w) :=
 \begin{cases}
  d_i(w) & \tau_{m(i)} < M(i)\\
  \tilde d_i(w) & \tau_{m(i)} \geq M(i).
\end{cases}
\end{equation}

\noindent As an example, we may take $\tau= 456667899$ and $w=534826179$, as in Figure \ref{tau}. Then $D_3^{(\tau)}(w)=\tilde d_3(w)=542836179 $ and $D_5^{(\tau)}(w)=d_5(w)=634825179$. 

We may generalize $\mg_n$ by defining $\mg^{(\tau)}_n$ as the $(n,n)$-signed colored graph with vertex set indexed by $S_n$, signature function given by inverse descents, and each edge set $E_i$ given by the nontrivial orbits of $D_i^{(\tau)}$.
Direct inspection shows that if $\tau=\tau(\bnu)$, then $D^{(\tau)}_i$ takes shifted content words of standard fillings of $\bnu$ to shifted content words of other standard fillings of $\bnu$. Thus, $D^{(\tau)}_i$ has a well-defined action on SYT$(\bnu)$ inherited from the action of $D^{(\tau)}_i$ on shifted content words. We may then define the following subgraph of $\mg_n^{(\tau)}$.

\begin{dfn}\label{dfn LLT graph}
\textup{Given some tuple of skew shapes $\bnu$, the \emph{LLT graph} $\ml_{\bnu} =(V, \sigma, E)$ is defined to be the $(n,n)$-signed colored graph with the following data:
\begin{enumerate}\itemsep.0in
\item $V$ 
$ = \{w\in S_n: w \textup{ is the shifted content word of some } {\bf T} \in \textup{SYT}({\bnu})\}$, 
\item The signature function $\sigma$ is given by the inverse descent sets of $w \in V$,
\item  The edge sets $E_i$ are defined by the nontrivial orbits of $D^{(\tau)}_i$ for all $1 < i < |\bnu|$, where $\tau=\tau(\bnu)$.
\end{enumerate}
}
\end{dfn}

\begin{example}\label{LLT example} \textup{ Consider $\bnu=((2),(2),(1),(1))$. A portion of the LLT graph $\ml_{\bnu}$ is presented in Figure \ref{fig LLT example}. Here, $\ml_{\bnu}$ is a subgraph of $\mg^{(\tau)}_6$ with $\tau = 566666.$ The entire connected component of the vertices in Figure \ref{fig LLT example} has 47 vertices. In the figure, the edge $\{312654, 412653\}$  is defined by the action of $d_3$ and $d_4$, while all other edges are defined by the action of $\tilde d_i$ for $1<i<6$.}

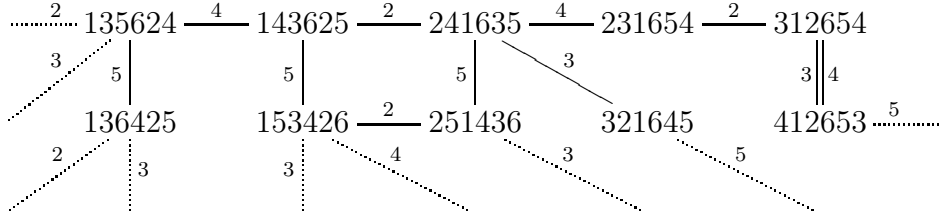
\begin{figure}[H]
\[
\xymatrix{
\ar@{.}[r]^{2\;\;\;\,}&{135624}\ar@{-}[r]^4\ar@{-}[d]_5& 143625\ar@{-}[d]_5\ar@{-}[r]^2 & 241635 \ar@{-}[r]^4\ar@{-}[d]_5 \ar@{-}[dr]^3 & 231654 \ar@{-}[r]^2& 312654 \ar@{=}[d]^4_3 & \\
\ar@{.}[ur]^3&136425 &153426 \ar@{-}[r]^2& 251436 & 321645 \ar@{.}[dr]^5  & 412653\ar@{.}[r]^{\;\;\;\,5} &\\
\ar@{.}[ur]^2&\ar@{.}[u]_3& \ar@{.}[u]^3 & \ar@{.}[ul]_4&  \ar@{.}[ul]_3 & & 
}
\]
\caption{A portion of $\ml_{\bnu}$ with signatures omitted. Here  $\bnu=((2), (2), (1), (1)).$}
\label{fig LLT example}
\end{figure}
\end{example}

\begin{remark}\label{LLT restrictions}
\textup{
\begin{enumerate} 
\parindent 0.25in
\item If $\bnu = (\Skew)$ is a 1-tuple, then each ${\bf T} \in \textup{SYT}(\bnu)$ corresponds to a unique $T\in\textup{SYT}({\Skew})$ in the obvious fashion. In particular, the shifted content word of ${\bf T}$ is equal to the content reading word of $T$. By Remark~\ref{skew via words},  $\mg_{\Skew} \cong \ml_\bnu$. In the case of straight shapes, $\bnu=(\lambda),$ Theorem~\ref{Skew} implies that $P\colon  S_n \rightarrow \textup{SYT}(n)$ induces an isomorphism from $\ml_\bnu$ to $\mg_\lambda$.
\item Given some $w$ in the vertex set of $\ml_\bnu$, we may readily describe the component of $w$ in a restriction of $\ml_\bnu$. Let $w$ be the shifted content word of ${\bf T}\in \textup{SYT}(\bnu)$, and let
 $\mc$ be the $(m,m)$-component of $w$ in $\ml_\bnu$. The isomorphism type of $\mc$ can be found by removing the cells of ${\bf T}$ containing values in $\{m+1, \ldots, n\}$,
 creating some ${\bf T^\prime}\in \textup{SYT}(\bnu^\prime)$ with shifted content word $w^\prime$. Then $\mc$ is isomorphic to the component of $w^\prime$ in $\ml_{\bnu^\prime}$. \\ \indent
Similarly, if $\mc$ is the component of $w$ in the $h$-upward restriction of $\ml_\bnu$, then the isomorphism type of $\mc$ can be found by removing the cells of ${\bf T}$ containing values in
$\{1, \ldots, h\}$ and subtracting $h$ from each of the remaining cells, creating some ${\bf T^\prime}\in \textup{SYT}(\bnu^\prime)$ with shifted content word $w^\prime$. Then $\mc$ is isomorphic to the component of $w^\prime$ in $\ml_{\bnu^\prime}$. \\ \indent
In both cases above, $\bnu^\prime$ is a tuple of skew shapes, so restrictions take components of LLT graphs to components of other LLT graphs. Furthermore, diam$(\bnu^\prime) \leq $ diam$(\bnu)$, since $\bnu^\prime$ is obtained by removing cells from $\bnu$.
\end{enumerate} }
\end{remark}
\leftskip 0.0in \parindent 0.25in

\noindent While LLT graphs do not necessarily satisfy Axiom 4 or Axiom 6, as can be seen in Example~\ref{LLT example}, they do satisfy a subset of the dual equivalence axioms. This is made precise in the following proposition and theorem.

\begin{prop}[\citep{Assaf} Prop.~4.6]\label{prop D}
Any LLT graph $\ml_\bnu$ obeys Axioms 1, 2, 3, and 5. Furthermore, the $\textup{inv}_\bnu$ statistic is constant on each connected component of $\ml_\bnu$.
\end{prop}

With Proposition \ref{prop D} in mind, it is natural to try to classify which LLT graphs satisfy all of the dual equivalence axioms. This classification is accomplished in the following theorem.

\begin{thm}\label{LLT graph}
The LLT graph $\ml_\bnu$ is a dual equivalence graph if and only if  \textup{diam}$(\bnu) \leq 3$.
\end{thm}

\noindent \emph{Proof:} First suppose diam$(\bnu) \leq 3$. By Proposition \ref{prop D}, $\ml_\bnu$ obeys Axioms 1, 2, 3, and 5. By Theorem~\ref{4plus}, we need only show that $\ml_\bnu$ satisfies Axiom $4^{+}$ to prove that $\ml_\bnu$ is a dual equivalence graph.

We begin by showing that $\ml_\bnu$ obeys Axiom $4$. Applying the same logic as Part~1 of Remark~\ref{relabeling 4}, we may use restrictions to only consider signed colored graphs of type $(5,5)$. By Part~2 of Remark~\ref{LLT restrictions}, each connected component of any such restriction is a component of $\ml_\bmu$ for some $\bmu$ with diam$(\bmu) \leq 3$. We may then reduce to the case where $|\bnu| = 5$ and diam$(\bnu)\leq 3$.  It thus suffices to consider $\mg_5^{(\tau)}$ for all $\tau=\tau(\bnu)$ with $\bnu$ satisfying these properties, in which case, it can be checked explicitly that all components of $\mg_5^{(\tau)}$ that are contained in some $\ml_\bnu$ with diam$(\bnu)\leq 3$ are dual equivalence graphs. This fact is verified via computer at 
\begin{center}
\textless http://www.math.washington.edu/$\sim$austinis/Proof\textunderscore LLTandDEG.sws\textgreater,
\end{center}
  cited as \citep{code}.

To demonstrate Axiom $4^+$, we may similarly inspect  $\mg^{(\tau)}_6$ for all $\tau$ that are derived from $\bnu$ with $|\bnu|=6$ and diam($\bnu)\leq 3$. It is possible to explicitly check that $\ml_\bnu$ is a dual equivalence graph by showing that the components of $\mg^{(\tau)}_6$ that satisfy Axiom 4 also satisfy the hypotheses of Corollary $\ref{computer check}$. This fact is also verified via computer at \citep{code}. Thus, $\ml_\bnu$ is a dual equivalence graph whenever diam$(\bnu)\leq 3$.

Now suppose diam$(\bnu)\geq 4$.  By the definition of diameter, there exist four cells of $\bnu$ whose shifted contents are distinct and differ by at most $k$. If we label these four cells $c_1, c_2, c_3$, and $c_4$ in shifted content reading order, then $c_1, c_2$, and $c_3$ must occur in different skew tableaux of $\bnu$. Similarly, $c_2, c_3$, and $c_4$ occur in different skew tableaux of $\bnu$. There is then some value of $i$ and some standard filling of $\bnu$ with values $i, i-1, i+2, i+1 $  or $ i+1, i+2, i-1, i $ in $c_1, c_2, c_3,$ and $c_4,$ respectively. Call the shifted content word of this standard filling $w$. Direct computation shows that $w$ and $D^{(\tau)}_i(w)=\tilde d_i(w)$ are contained in distinct $(i+1)$-edges. That is, the restriction of $\ml_\bnu$ to the edge colored graph $(V, E_i\cup E_{i+1})$ contains a component with at least three distinct edges. Therefore, $\ml_\bnu$ violates Axiom 4 and is not a dual equivalence graph.~\hfill~$\square$


\begin{remark}
\textup{
By applying \cite[Prop.~5.3]{Assaf08}, it can be proven that $\ml_{\bnu}$ is a dual equivalence graph whenever $\bnu$ is a pair of ribbons. Both the previous theorem and the following lemma can be viewed as extensions of this fact.
}
\end{remark}

While the previous theorem describes which LLT graphs are dual equivalence graphs, there are other cases in which only specific components of an LLT graph are dual equivalence graphs. Sometimes, we may even give an explicit isomorphism from a component of an LLT graph to a standard dual equivalence graph. 
The next three results describe cases where such an explicit isomorphism exists.

\begin{lem}\label{gap}
Let $\mc=(V,\sigma, E)$ be a component of $\mg_n^{(\tau)}$ such that $D^{(\tau)}_i$ never acts on $w\in V$ nontrivially via $\tilde d_i$ unless $i-1, i,$ and $i+1$ have adjacent indices in $w$. Then $\mc$ is a dual equivalence graph, and $P \hspace{-.04in}:\hspace{-.04in}V \rightarrow \textup{SYT}(n)$ induces an isomorphism from $\mc$ to some standard dual equivalence graph $\mg_\lambda$.
\end{lem}

\noindent $Proof:$ By Proposition \ref{prop D}, $\mc$ satisfies Axioms 1, 2, 3, and 5. By Theorem~\ref{4plus}, we need only demonstrate that $\mc$ satisfies Axiom $4^+$ to show that $\mc$ is a dual equivalence graph. As in the proof of Theorem~\ref{LLT graph}, we may reduce to the case where $\mc$ has type (6,6).  The requirement on the action of $D^{(\tau)}_i$ allows us to further reduce to the case where $\tau$ satisfies $\tau_i \leq i+2$ for all $1\leq i \leq n$. In order to show that $\mc$ is a dual equivalence graph, it is thus sufficient to check that $\mg_6^{(\tau)}$ is a dual equivalence graph for all $\tau$ of length $6$ with $\tau_i \leq i+2$. This fact is verified via computer at \citep{code}.

To demonstrate that $P$ induces a morphism, we show that $\mc$ satisfies the hypotheses of Theorem~\ref{KnuthMorph}. That is, we need to show that $E_i(w)$ is in the Knuth class of $d_i(w)$ for all $w\in V$ and all $1<i<n$. When $D^{(\tau)}_i$ acts via $d_i$, this edge requirement is clearly satisfied. Now notice that if $D^{(\tau)}_i$ acts via $\tilde d_i$, then the restriction on indices implies that $\tilde d_i=K_j \circ d_i$ for some $j$. Hence, $E_i(w)$ is in the same Knuth class as $d_i(w)$. 
Thus, $P$ induces a morphism between connected dual equivalence graphs. This morphism must be an isomorphism by Corollary~\ref{terminal}.~\hfill~$\square$
 
\vspace{.1in}

Lemma~\ref{gap} may be readily applied to LLT graphs because every LLT graph is a subgraph of some $\mg_n^{(\tau)}$. While this lemma does give a particularly nice isomorphism, it does not apply in the more general case of LLT graphs indexed by tuples of shapes with diameter at most 3, as in Theorem~\ref{LLT Decomposition}. Instead, we will first need to identify specific components of LLT graphs satisfying the hypotheses of Lemma~\ref{gap}.

\begin{lem}\label{two row lemma}
Let $\bnu$ be a tuple of skew shapes such that $|\bnu|=n$ and $\textup{diam}(\bnu)\leq 3$. Let $\mc=(V, \sigma, E)$ be a connected component of $\ml_{\bnu}$, let $\tau=\tau(\bnu)$, and let $v$ be any vertex in $\mc$. If $\textup{sh}(v)$ has strictly less than three rows, then $D^{(\tau)}_i$ never acts on any vertex $w\in V$ nontrivially via $\tilde d_i$ unless $i-1, i,$ and $i+1$ have adjacent indices in $w$.
\end{lem}

\noindent $Proof:$ Let $v$ be the shifted content word of ${\bf T}=(T^{(0)}, \ldots, T^{(k-1)}) \in \textup{SYT}(\bnu)$. We claim that it is enough to show that $D^{(\tau)}_i$ never acts on $v$ nontrivially via $\tilde d_i$ unless $i-1, i,$ and $i+1$ have adjacent indices in $v$. In this case, $D^{(\tau)}_i$ acts on $v$ as the identity, $d_i$, or $K_j\circ d_i$ for some $j$, as was mentioned in the proof of Lemma~\ref{gap}. Thus sh($D^{(\tau)}_i(v))= $ sh$(v)$. In particular,  $D^{(\tau)}_i(v)$ satisfies the same hypotheses as $v$ and recursively, so does every vertex of $\mc$. 

We may recharacterize the condition that sh$(v)$ has strictly less than three rows in terms of the values in $v$. If $P(v)$ has strictly less than three rows, then  $v$ has no decreasing subword of length 3, as noted in Section \ref{The RSK Correspondence}. Suppose $D^{(\tau)}_i(v)=\tilde d_i(v) \neq v$ for some fixed choice of $i$, and let $a$, $b$, and $c$ be the values $i-1, i,$ and $i+1$ given in the order they appear in $v$. Because $D^{(\tau)}_i$ acts nontrivially on $v$, either $a > b$ or $b > c$. Hence, there cannot be any value strictly greater than $i+1$ that occurs before $a$ in $v$ or any value strictly less than $i-1$ that occurs after $c$ in $v$.

We can now show that $a, b$, and $c$ occur consecutively in $v$ by considering shifted contents. Let $\tilde c(x)$ denote the shifted content of the cell of ${\bf T}$ containing the value $x$. Suppose that there is a value $x$ that appears after $a$ in $v$ with $\tilde c(a)= \tilde c(x)$. That is, $x$ occurs above $a$ on the same diagonal of some skew tableau $T^{(j)}$ of ${\bf T}$. Then there is a value $y$ that occurs directly north of $a$ and directly west of $x$ in $T^{(j)}$. It follows that $y> i+1 \geq a$, and $y$ occurs before $a$ in $v$, a contradiction. Similarly, there cannot be a value $x$ that occurs before $c$ in $v$ with $\tilde c(c)= \tilde c(x)$. In particular, $\tilde c(a) < \tilde c(b) < \tilde c(c)$.

Because diam$(\bnu) \leq 3$, the cells of $\bnu$ can occupy at most one shifted content strictly between $\tilde c(a)$ and $\tilde c(c)$, so it suffices to show that there is exactly one cell in $\bnu$ with shifted content $\tilde c(b)$. Otherwise, there must be a cell on the same diagonal as $b$ in some skew tableau of ${\bf T}$. There is then a value $z$ in a cell either directly north of $b$ or directly south of $b$ with $\tilde c(z)=\tilde c (b)-k$ or $\tilde c(z)=\tilde c (b)+k$, respectively. Also, $z > b$ or $z < b$, respectively.
Because $D_i^{(\tau)}(v)=\tilde d_i(v)\neq v$, $\tilde c (c)- \tilde c(a)\leq k$. Using the fact that $\tilde c(a) < \tilde c(b) < \tilde c(c)$, it follows that $\tilde c(b)-k < \tilde c(a)$ and $\tilde c(c) -\tilde c(b) < k$. Thus, $z$ is greater than $i+1$ and occurs before $a$ in $v$, or $z$ is less than $i-1$ and occurs after $c$ in $v$, forcing a contradiction.
Therefore, $a, b,$ and $c$ occur consecutively in $v$.~\hfill~$\square$
\vspace{.1in}

\begin{cor}\label{two row case}
Let $\bnu$ be a tuple of skew shapes such that $|\bnu|=n$ and $\textup{diam}(\bnu)\leq 3$. Let $\mc=(V, \sigma, E)$ be a connected component of $\ml_{\bnu}$, and let $w$ be any vertex in $\mc$. If $\lambda=\textup{sh}(w)$ has strictly less than three rows, then $P\colon  V \rightarrow SYT(n)$ induces an isomorphism from $\mc$ to the standard dual equivalence graph $\mg_\lambda$.\end{cor}

\noindent $Proof:$ The corollary follows immediately from Lemma~\ref{gap} and Lemma~\ref{two row lemma}.

\vspace{.1in}

We now change our focus from finding isomorphism types to finding a set of vertices to represent the components of an LLT graph. The following lemma is crucial to the proof of Theorem~\ref{LLT Decomposition}.

\begin{lem}\label{Luv}
Let $\bmu$ and $\bnu$ be tuples of skew shapes such that $\textup{diam}(\bmu), \textup{diam}(\bnu)\leq 3$ and $|\bmu|=|\bnu|=n$. Let $\mc$ and $\md$ be connected components of $\ml_\bmu$ and $\ml_\bnu$, respectively, and let $\phi\colon \mc\rightarrow \md$ be an isomorphism. If $w$ is a vertex in $\mc$ and $\lambda \vdash n$, then $w\in \textup{SYam}(\lambda)$ if and only if $\phi(w) \in \textup{SYam}(\lambda)$. 
\end{lem}

\noindent $Proof:$   Suppose, for the sake of contradiction, that $v=\phi(w) \in \textup{SYam}(\lambda)$ and $w \notin \textup{SYam}(\lambda)$. We begin by setting some definitions. Let $P(w)=T$. In particular, $P(\phi(w)) = U_\lambda \neq T$. 
Consider the lowest value that does not occur in the same cell of $T$ and $U_\lambda$. By signature considerations, this value must be in a lower row of $T$ than in $U_\lambda$. Let $m$ be the smallest number such that some entry of the $m^{th}$ row of $U_\lambda$ occurs in a lower row of $T$ (see Figure \ref{2row2} for an example).  
Now define
\begin{equation}
p = \sum_{j=1}^{m-2} \lambda_j, \hspace{.22in}
q = \lambda_{m-1}+\lambda_m, \hspace{.22in} 
S_i = p+i,   \hspace{.22in}
\textup{and} \hspace{.22in}
S =\{ S_1, S_2, \ldots, S_q\}.
\end{equation}
That is, $S$ is the set of values in rows 
$m-1$ and $m$ of $U_\lambda$. 
Notice that $m >1$, so $S$ is nonempty. 
Let $w_S$ and $v_S$ be the subwords of $w$ and $v$, respectively, with values in $S$. 

We now consider the longest increasing subwords of $w_S$ and $v_S$. Because $v$ is Knuth equivalent to the row reading word of $U_\lambda$, we may restrict the row reading word of $U_\lambda$ to values in $S$ in order to find the longest increasing subword of $v_S$, 
as described in Section \ref{The RSK Correspondence}. Specifically, the longest increasing subword of $v_S$ has length $\lambda_{m-1}$. Similarly, we may restrict the row reading word of $T$ to find the longest increasing subword of $w_S$.
 By the definitions of $m$ and $S$, the values $S_1$ through $S_{\lambda_{m-1}}$ occur in row $m-1$ of $T$, while the value $S_q$ occurs no higher than row $m-1$ of $T$. Hence, the row reading word of $T$ has 
$S_1 S_2 \ldots S_{\lambda_{m-1}} S_{q}$ as an increasing subword (see Figure \ref{2row2} for an example).
In particular, if $l$ is the length of the longest increasing subword of $w_S$, then $l > \lambda_{m-1}$. We will use these facts about $w_S$ and $v_S$ to create a contradiction to the assumption that $T \neq U_\lambda$. 

\begin{figure}[H]
   \begin{center}
        \ytableausetup{aligntableaux=center}
	$T= \ytableaushort{{\bf{7}},{\bf{4}}{\bf{5}}{\bf{6}},123{\bf{8}}{\bf{9}}}  \hspace{.8in}
	U_{(3,3,3)}=\ytableaushort{{\bf{7}}{\bf{8}}{\bf{9}},{\bf{4}}{\bf{5}}{\bf{6}},123}$
		\end{center}
	\caption{An example where $T$ has the same signature as $U_{(3,3,3)}$. Here, $m=3, p=3,$  \hspace{.1in} $q=6, S=\{4, 5, 6, 7, 8, 9\}$, $l = 5$, and  the subword $S_1S_2\ldots S_{\lambda_{m-1}}S_q= 4569$.}
	\label{2row2}
\end{figure}

Consider the  $(q, q)$-component of $w$ in the $p$-upward restriction of $\mc$. Call this component $\mc^\prime$. We proceed by finding the isomorphism type of $\mc^\prime$ in two different ways.
First, 
we may find the isomorphism type of $\mc^\prime$ in the manner described by Part~2 of Remark~\ref{LLT restrictions}. Specifically, consider ${\bf T} \in \textup{SYT}(\bmu)$ such that $w$ is the shifted content word of ${\bf T}$. 
Remove the cells of ${\bf T} $ not containing values in $S$ and then subtract $p$ from each of the remaining cells, creating some ${\bf T}^\prime$ of shape $\bmu^\prime$ with shifted content word $w^\prime$.  Then $\mc^\prime$ is isomorphic to the component of $\ml_{\bmu^\prime}$ containing $w^\prime$. Here, $w^\prime$ can also be found by subtracting $p$ from each entry in $w_S$. 

To find the isomorphism type of the component of $w^\prime$ in $\ml_{\bmu^\prime}$, we will apply Corollary~\ref{two row case}.
The signature of $w^\prime$ is equal to the restriction of the signature of $w$ to the coordinates $S_1$ through $S_{q-1}$. The signature of $w$ is equal to the signature of $v$, which is equal to the signature of $U_\lambda$. Thus, the signature of $w^\prime$ is equal to the restriction of the signature of $U_\lambda$ to the coordinates $S_1$ through $S_{q -1}$. In particular, the signature of $w^\prime$ has exactly one $-1$, implying that $P(w^\prime)$ has exactly two rows. Applying Corollary~\ref{two row case}, $P$ induces an isomorphism from $\mc^\prime$ to the standard dual equivalence graph $\mg_{\textup{sh}(w^\prime)}$.  The length of the first row in sh$(w^\prime)$ is the length of the longest increasing subword in $w^\prime$, as noted in Section \ref{The RSK Correspondence}. The length of this subword is equal to the length of the longest increasing subword of $w_S$, which we know to have length $l>\lambda_{m-1}$. Since sh$(w^\prime)$ has exactly two rows, it follows that sh($w^\prime)=(l, q-l)$. Therefore, $\mc^\prime \cong \mg_{(l, q-l)}$. 
 
Alternatively, we may find the isomorphism type of $\mc^\prime$ by applying the same restrictions to $\md$. Let $\md^\prime$ be the $(q, q)$-component of $v$ in the $p$-upward restriction of $\md$. In particular, $\mc^\prime \cong \md^\prime$. As before, $\md^\prime$ is isomorphic to the component of some $\ml_{\bnu^\prime}$ containing  the vertex $v^\prime$, where $v^\prime$ is obtained by subtracting $p$ from each of the values in $v_S$. Also as before, sh$(v^\prime)$ has two rows. Therefore Corollary~\ref{two row case} guarantees that $D^\prime \cong \mg_{\textup{sh}(v^\prime)}$. The length of the longest increasing subword of $v^\prime$ is equal to the length of the longest increasing subword of $v_S$, which we know to have length $\lambda_{m-1}$. Thus, sh$(v^\prime) = (\lambda_{m-1}, \lambda_{m})$. Since $l > \lambda_{m-1}$, $\mc^\prime \cong \mg_{(l, q-l)} \ncong \mg_{(\lambda_{m-1}, \lambda_{m})} \cong \md^\prime$, a contradiction. Therefore, $w \in \textup{SYam}(\lambda)$ whenever $v=\phi(w)\in \textup{SYam}(\lambda)$.

We still need to prove that $\phi(w) \in \textup{SYam}(\lambda)$ whenever $w \in \textup{SYam}(\lambda)$. This follows via symmetry by considering $w$ as the image of $\phi(w)$ under the isomorphism $\phi^{-1}$. Thus, $w \in \textup{SYam}(\lambda)$ if and only if $\phi(w) \in \textup{SYam}(\lambda)$.~\hfill~$\square$

\vspace{.1in}
\noindent 
\textbf{Proof of Theorem~\ref{LLT Decomposition}:}
~\\
\indent We begin by reducing Theorem~\ref{LLT Decomposition} to a statement about signed colored graphs.
Let $V_1, V_2, \ldots, V_m$ be the vertex sets of the connected components of $\ml_\bnu = (V, \sigma, E)$. Applying (\ref{dfn LLT poly}) and Definition~\ref{dfn LLT graph},
\begin{align}\label{llt sum}
\llt =
 \sum_{{\bf T} \in \textup{SYT}(\bnu)} q^{\textup{inv}_k({\bf T})} F_{\sigma({\bf T})}(X) = 
\sum_{v\in V} q^{\textup{inv}_\bnu(v)} F_{\sigma(v)}(X) = 
\sum_{j=1}^m \sum_{v\in V_j} q^{\textup{inv}_\bnu(v)} F_{\sigma(v)}(X), 
\end{align}
where it should be noted that the signature function $\sigma$ changes between the first and second sums. Specifically, $\sigma$ changes from the signature function on standard fillings of $\bnu$ to the signature function given by the inverse descents of the permutations in the vertex set of $\ml$. The second equality then follows because the signature of a standard filling of $\bnu$ is defined via the signature of its row reading word.

To further simplify (\ref{llt sum}), we turn our attention to the individual components of $\ml_\bnu$ when diam$(\bnu)\leq3$. Let $\mc=(V_j, \sigma, E)$ be a connected component of $\ml_\bnu$, and choose any fixed $v_j\in V_j$. By Proposition \ref{prop D}, the inv$_\bnu$ statistic is constant on $V_j$. By Theorems \ref{main} and \ref{LLT graph}, $\mc$ is isomorphic to $\mg_\lambda$ for some $\lambda \vdash |\bnu|$. From (\ref{schurs}), it follows that
\begin{equation}
\sum_{v\in V_j} q^{\textup{inv}_\bnu(v)} F_{\sigma(v)}(X) = q^{\textup{inv}_\bnu(v_j)} s_\lambda.
\end{equation}
To prove Theorem \ref{LLT Decomposition}, we need only guarantee that there is a unique choice of $v_j\in V_j$ such that $v_j$ is a standardized Yamanouchi word and, moreover, that this choice of $v_j$ is in $\textup{SYam}(\lambda)$.

We will find such a $v_j$ explicitly by applying Lemma \ref{Luv}. 
Consider the LLT graph $\ml_\blambda$, where $\blambda$ is the 1-tuple $(\lambda)$. It follows from Part~2 of Remark~\ref{LLT restrictions} that diam$(\blambda) \leq 2$. By Part~1 of Remark~\ref{LLT restrictions}, $P$ induces an isomorphism from $\ml_\blambda$ to $\mg_\lambda$. In particular, the vertex set of $\ml_\blambda$ has exactly one vertex $w$ such that $P(w) = U_\lambda$. Thus, $w\in \textup{SYam}(\lambda)$, and $w$ is the only vertex of $\ml_\blambda$ that is a standardized Yamanouchi word. 
Because $\mc$ is isomorphic to $\mg_\lambda$, and $\mg_\lambda$ is isomorphic to $\ml_\blambda$, there exists an isomorphism $\phi\colon  \ml_\blambda \rightarrow \mc$. Applying Lemma~\ref{Luv}, $v_j = \phi(w)$ is the unique standardized Yamanouchi word in $V_j$.~\hfill~$\square$

\section{Conclusion}

There are a number of persistent open questions involving dual equivalence graphs. We conclude by mentioning some of these questions as possibilities for further research. 

Can dual equivalence graphs be used to give nice Schur expansions for symmetric functions other than LLT and Macdonald polynomials? Many functions have known expansions in terms of fundamental quasisymmetric functions but are lacking nice descriptions for their Schur expansions. As was the case with LLT polynomials, it may be possible to apply dual equivalence graphs to find Schur expansions for such functions. Examples include plethysms of Schur functions and k-Schur functions. The former is described in \citep{LW}, while the latter has already seen some progress in \citep{AB}. Though currently unproven, the shuffle conjecture provides another example. This conjecture describes the composition of the nabla operator with an elementary symmetric function in terms of two statistics and fundamental quasisymmetric functions. This sum, in turn, is equal to a sum of LLT polynomials with an additional statistic. It may then be possible to use the results of Section \ref{LLTs} as an aid to understanding the combinatorics of the shuffle conjecture.  For the original statement of the shuffle conjecture, along with proofs of the facts mentioned above, see \citep{HHLRU}.

Is there a good description of when a component of an LLT graph is a dual equivalence graph?
While Theorem~\ref{LLT graph} classifies when an LLT graph is a dual equivalence graph, it makes no claims about specific components of LLT graphs in the diam$(\bnu)\geq 4$ case.  
Similarly, there is currently no good description for when a specific component of $\mg_n^{(\tau)}$ is a dual equivalence graph.

Is there an axiomatic description for when a signed colored graph admits a morphism onto a standard dual equivalence graph?
While Theorem~\ref{KnuthMorph} gives one criterion for when a signed colored graph admits a morphism onto a standard dual equivalence graph, a more axiomatic description would be desirable. Such a description would necessarily provide sufficient conditions for when a signed colored graph corresponds to a positive integer multiple of a Schur function. 

Can the axiomatization of dual equivalence graphs be generalized to an axiomatization of the family of signed colored graphs that correspond to Schur positive functions? Such an axiomatization would necessarily be less strict than the one given for dual equivalence graphs and would thus be satisfied by a larger set of signed colored graphs.
%
%
 An axiomatization of this family of graphs could provide expanded methods for proving the Schur positivity of a variety of symmetric functions.


\bibliographystyle{abbrvnat}

\begin{thebibliography}{30}
\providecommand{\natexlab}[1]{#1}
\providecommand{\url}[1]{\texttt{#1}}
\expandafter\ifx\csname urlstyle\endcsname\relax
  \providecommand{\doi}[1]{doi: #1}\else
  \providecommand{\doi}{doi: \begingroup \urlstyle{rm}\Url}\fi

\bibitem[Assaf(2007)]{Assaf07}
S.~H. Assaf.
\newblock \emph{Dual equivalence graphs, ribbon tableaux and {M}acdonald
  Polynomials}.
\newblock PhD thesis, University of California at Berkeley, 2007.

\bibitem[Assaf(2008/09)]{Assaf08}
S.~H. Assaf.
\newblock A generalized major index statistic.
\newblock \emph{S\'{e}m. Lothar. Combin}, 60, 2008/09.

\bibitem[Assaf(2011)]{Assaf}
S.~H. Assaf.
\newblock Dual equivalence graphs {I}: A combinatorial proof of {LLT} and
  {M}acdonald polynomials.
\newblock arXiv:1005.3759v2 [math.CO], 2011.

\bibitem[Assaf and Billey(2012)]{AB}
S.~H. Assaf and S.~Billey.
\newblock Affine dual equivalence and k-{S}chur functions.
\newblock \emph{J. Combinatorics}, 3\penalty0 (3):\penalty0 343--399, 2012.

\bibitem[Assaf et~al.(2012)Assaf, Bergeron, and Sottile]{ABS}
S.~H. Assaf, N.~Bergeron, and F.~Sottile.
\newblock Multiplying {S}chubert polynomials by {S}chur polynomials.
\newblock The 5th Mathematical Society of Japan Seasonal Institute, 2012.
\newblock (manuscript).

\bibitem[Carr\'{e} and Leclerc(1995)]{CL}
C.~Carr\'{e} and B.~Leclerc.
\newblock Splitting the square of a {S}chur function into its symmetric and
  antisymmetric parts.
\newblock \emph{J. Algebraic. Combin.}, 4\penalty0 (3), 1995.

\bibitem[Fishel(1995)]{Fishel}
S.~Fishel.
\newblock Statistics for special q,t-{K}ostka polynomials.
\newblock \emph{Proc. Amer. Math. Soc}, 123\penalty0 (10):\penalty0 2961--2969,
  1995.

\bibitem[Fulton(1997)]{Fulton}
W.~Fulton.
\newblock \emph{Young Tableaux, with Applications to Representation Theory and
  Geometry}.
\newblock London Mathematical Society Student Texts 35, Cambridge University
  Press, Cambridge, 1997.

\bibitem[Gessel(1984)]{Gessel}
I.~M. Gessel.
\newblock Multipartite {P}-partitions and inner products of skew {S}chur
  functions.
\newblock \emph{Contemp. Math}, 34:\penalty0 289--317, 1984.

\bibitem[Godsil and Royle(2004)]{GR}
C.~Godsil and G.~Royle.
\newblock \emph{Algebraic Graph Theory}, volume 207 of \emph{Graduate Texts in
  Mathematics}.
\newblock Springer-Verlag, New York, 2004.

\bibitem[Haglund(2004)]{HaglundCombinatorialModel}
J.~Haglund.
\newblock A combinatorial model for the {M}acdonald polynomials.
\newblock \emph{Proc. Natl. Acad. Sci. USA}, 101\penalty0 (46):\penalty0
  16127--16131, 2004.

\bibitem[Haglund(2008)]{Haglund}
J.~Haglund.
\newblock \emph{The {$q$},{$t$}-{C}atalan Numbers and the Space of Diagonal
  Harmonics}, volume~41 of \emph{University Lecture Series}.
\newblock American Mathematical Society, Providence, RI, 2008.

\bibitem[Haglund et~al.(2005{\natexlab{a}})Haglund, Haiman, and Loehr]{HHL}
J.~Haglund, M.~Haiman, and N.~Loehr.
\newblock A combinatorial formula for {M}acdonald polynomials.
\newblock \emph{J. Amer. Math. Soc}, 18\penalty0 (3), 2005{\natexlab{a}}.

\bibitem[Haglund et~al.(2005{\natexlab{b}})Haglund, Haiman, Loehr, Remmel, and
  Ulyanov]{HHLRU}
J.~Haglund, M.~Haiman, N.~Loehr, J.~B. Remmel, and A.~Ulyanov.
\newblock A combinatorial formula for the character of the diagonal
  coinvariants.
\newblock \emph{Duke Math. J}, 126\penalty0 (2):\penalty0 195--232,
  2005{\natexlab{b}}.

\bibitem[Haiman(2001)]{Hai01}
M.~Haiman.
\newblock Hilbert schemes, polygraphs and the {M}acdonald positivity
  conjecture.
\newblock \emph{J. Amer. Math. Soc}, 14\penalty0 (4):\penalty0 941--1006, 2001.

\bibitem[Haiman(1992)]{Haiman}
M.~D. Haiman.
\newblock Dual equivalence with applications, including a conjecture of
  {P}roctor.
\newblock \emph{Discrete Math}, 99\penalty0 (1):\penalty0 79--113, 1992.

\bibitem[Knuth(1970)]{Knuth}
D.~E. Knuth.
\newblock Permutations, matrices, and generalized {Y}oung tableaux.
\newblock \emph{Pacific J. Math.}, 34:\penalty0 709--727, 1970.

\bibitem[Lapointe and Morse(2003)]{LM}
L.~Lapointe and J.~Morse.
\newblock Tableaux statistics for two part {M}acdonald polynomials.
\newblock \emph{Algebraic combinatorics and quantum groups}, pages 61--84,
  2003.

\bibitem[Lascoux et~al.(1997)Lascoux, Leclerc, and Thibon]{LLT}
A.~Lascoux, B.~Leclerc, and J.-Y. Thibon.
\newblock Ribbon tableaux, {H}all-{L}ittlewood functions, quantum affine
  algebras, and unipotent varieties.
\newblock \emph{J. Math. Phys}, 38\penalty0 (2):\penalty0 1041--1068, 1997.

\bibitem[Loehr and Warrington(2012)]{LW}
N.~A. Loehr and G.~S. Warrington.
\newblock Quasisymmetric expansion of {S}chur function plethysms.
\newblock \emph{Proc. Amer. Math. Soc}, 140\penalty0 (4):\penalty0 1159--1171,
  2012.

\bibitem[Macdonald(1988)]{Mac}
I.~G. Macdonald.
\newblock A new class of symmetric functions.
\newblock \emph{Actes du $20^e$ S\'{e}minaire Lotharingien}, 372:\penalty0
  131--171, 1988.

\bibitem[Macdonald(1995)]{Macdonald}
I.~G. Macdonald.
\newblock \emph{Symmetric functions and {H}all polynomials}.
\newblock Oxford Mathematical Monographs. The Clarendon Press Oxford University
  Press, New York, second edition, 1995.

\bibitem[Roberts(2012)]{code}
A.~Roberts.
\newblock Computer proofs: {LLT} polynomials and dual equivalence graphs.
\newblock $<$http://www.math.washington.edu/$\sim$austinis/
  Proof\textunderscore LLTandDEG.sws$>$, 2012.

\bibitem[Sagan(2001)]{Sagan}
B.~E. Sagan.
\newblock \emph{The Symmetric Group: Representations, Combinatorial Algorithms,
  and Symmetric Functions}, volume 203 of \emph{Graduate Texts in Mathematics}.
\newblock Springer-Verlag, New York, second edition, 2001.

\bibitem[Stanley(1999)]{Stanley2}
R.~P. Stanley.
\newblock \emph{Enumerative Combinatorics. {V}ol. 2}, volume~62 of
  \emph{Cambridge Studies in Advanced Mathematics}.
\newblock Cambridge University Press, Cambridge, 1999.

\bibitem[Stanley(2012{\natexlab{a}})]{Catalan}
R.~P. Stanley.
\newblock \emph{Catalan Addendum}.
\newblock $<$www-math.mit.edu/ rstan/ec/catadd.pdf$>$, 2012{\natexlab{a}}.

\bibitem[Stanley(2012{\natexlab{b}})]{Stanley}
R.~P. Stanley.
\newblock \emph{Enumerative Combinatorices. {V}ol. 1}, volume~49 of
  \emph{Cambridge Studies in Advanced Mathematics}.
\newblock Cambridge University Press, Cambridge, second edition,
  2012{\natexlab{b}}.

\bibitem[van Leeuwen(2000)]{vanL}
M.~A.~A. van Leeuwen.
\newblock Some bijective correspondences involving domino tableaux.
\newblock \emph{Electron. J. Combin}, 7:\penalty0 Research Paper 35, 25 pp
  (electronic), 2000.

\bibitem[Zabrocki(1998)]{Zab}
M.~Zabrocki.
\newblock A {M}acdonald vertex operator and standard tableaux statistics for
  the two-column {$(q,t)$}-{K}ostka coefficients.
\newblock \emph{Electron. J. Combin.}, 5:\penalty0 Research Paper 45, 46, 1998.

\bibitem[Zabrocki(1999)]{ZabrockiSpecialCase}
M.~Zabrocki.
\newblock Positivity for special cases of {$(q,t)$}-{K}ostka coefficients and
  standard tableaux statistics.
\newblock \emph{Electron. J. Combin.}, 6:\penalty0 Research Paper 41, 36, 1999.

\end{thebibliography}

\end{document}